  %
  %

  \magnification 1200


  \newcount\fontset
  \fontset=1
  \def\dualfont#1#2#3{\font#1=\ifnum\fontset=1 #2\else#3\fi}

  \dualfont\bbfive{bbm5}{cmbx5}
  \dualfont\bbseven{bbm7}{cmbx7}
  \dualfont\bbten{bbm10}{cmbx10}
  \font \eightbf = cmbx8
  \font \eighti = cmmi8 \skewchar \eighti = '177
  \font \eightit = cmti8
  \font \eightrm = cmr8
  \font \eightsl = cmsl8
  \font \eightsy = cmsy8 \skewchar \eightsy = '60
  \font \eighttt = cmtt8 \hyphenchar\eighttt = -1
  \font \msbm = msbm10
  \font \sixbf = cmbx6
  \font \sixi = cmmi6 \skewchar \sixi = '177
  \font \sixrm = cmr6
  \font \sixsy = cmsy6 \skewchar \sixsy = '60
  \font \tensc = cmcsc10
  
  \font \titlefont = cmr7 scaled \magstep4
  \scriptfont \bffam = \bbseven
  \scriptscriptfont \bffam = \bbfive
  \textfont \bffam = \bbten

  \font\rs=rsfs10 
  \font\rssmall=rsfs8

  \newskip \ttglue

  \def \eightpoint {\def \rm {\fam0 \eightrm }%
  \textfont0 = \eightrm
  \scriptfont0 = \sixrm \scriptscriptfont0 = \fiverm
  \textfont1 = \eighti
  \scriptfont1 = \sixi \scriptscriptfont1 = \fivei
  \textfont2 = \eightsy
  \scriptfont2 = \sixsy \scriptscriptfont2 = \fivesy
  \textfont3 = \tenex
  \scriptfont3 = \tenex \scriptscriptfont3 = \tenex
  \def \it {\fam \itfam \eightit }%
  \textfont \itfam = \eightit
  \def \sl {\fam \slfam \eightsl }%
  \textfont \slfam = \eightsl
  \def \bf {\fam \bffam \eightbf }%
  \textfont \bffam = \eightbf
  \scriptfont \bffam = \sixbf
  \scriptscriptfont \bffam = \fivebf
  \def \tt {\fam \ttfam \eighttt }%
  \textfont \ttfam = \eighttt
  \tt \ttglue = .5em plus.25em minus.15em
  \normalbaselineskip = 9pt
  \def \MF {{\manual opqr}\-{\manual stuq}}%
  \let \sc = \sixrm
  \let \big = \eightbig
  \setbox \strutbox = \hbox {\vrule height7pt depth2pt width0pt}%
  \normalbaselines \rm }


 \def \Headlines #1#2{\nopagenumbers
    \advance \voffset by 2\baselineskip
    \advance \vsize by -\voffset
    \headline {\ifnum \pageno = 1 \hfil
    \else \ifodd \pageno \tensc \hfil \lcase {#1} \hfil \folio
    \else \tensc \folio \hfil \lcase {#2} \hfil
    \fi \fi }}
  \def \Date #1 {\footnote {}{\eightit Date: #1.}}


  \def \lcase #1{\edef \auxvar {\lowercase {#1}}\auxvar }

  \def \goodbreak {\vskip0pt plus.1\vsize \penalty -250 \vskip0pt
plus-.1\vsize }

  \newcount \secno \secno = 0
  \newcount \stno

  \def \seqnumbering {\global \advance \stno by 1
    \number \secno .\number \stno }

  \def\section #1{\global\def\SectionName{#1}\stno = 0 \global
\advance \secno by 1 \bigskip \bigskip \goodbreak \noindent {\bf
\number \secno .\enspace #1.}\medskip \noindent \ignorespaces}

  \long \def \sysstate #1#2#3{\medbreak \noindent {\bf \seqnumbering
.\enspace #1.\enspace }{#2#3\vskip 0pt}\medbreak }
  \def \state #1 #2\par {\sysstate {#1}{\sl }{#2}}
  \def \definition #1\par {\sysstate {Definition}{\rm }{#1}}
  \def \remark #1\par {\sysstate {Remark}{\rm }{#1}}


  \def \proof {\medbreak \noindent {\it Proof.\enspace }}
  \def \proofend {\ifmmode \eqno \square \else \hfill \square
\looseness = -1 \medbreak \fi }

  \def \$#1{#1 $$$$ #1}
  \def\=#1{\buildrel #1 \over =}

  \def\Item #1{\smallskip \item {#1}}
  \newcount \zitemno \zitemno = 0
  \def\izitem {\zitemno = 0}
  \def\zitem {\global \advance \zitemno by 1 \Item {{\rm(\romannumeral
\zitemno)}}}

  \newcount \footno \footno = 1
  \newcount \halffootno \footno = 1
  \def\footcntr {\global \advance \footno by 1
  \halffootno =\footno
  \divide \halffootno by 2
  $^{\number\halffootno}$}
  \def\fn#1{\footnote{\footcntr}{\eightpoint#1}}


  \def \N {{\bf N}}
  \def \<{\left \langle }   \def \<{\langle }
  \def \>{\right \rangle }  \def \>{\rangle }
  \def \curly#1{\hbox{\rs #1\/}}
  \def \ds{\displaystyle}
  \def \and {\hbox {,\quad and \quad }}
  \def \calcat #1{\,{\vrule height8pt depth4pt}_{\,#1}}
  
  \def \for #1{,\quad \forall\,#1}
  \def \square {\hbox {$\sqcap \!\!\!\!\sqcup $}}
  \def \crossproduct {{\hbox {\msbm o}}}
  \def \stress #1{{\it #1}\/}
  \def \inv {^{-1}}
  \def \* {\otimes}



  \def \ifn #1{\expandafter \ifx \csname #1\endcsname \relax }
  \def \track #1#2#3{\ifn{#1}\else { \tt #2\string #3 }\fi}
  \def \cite #1{{\rm [\bf #1\track{showref}{\#}{#1}\rm]}}
  \def \scite #1#2{{\rm [\bf #1\track{showref}{\#}{#1}{\rm \hskip 0.7pt:\hskip 2pt #2}\rm]}}
  \def \label #1{\global \edef #1{\number \secno \ifnum \number \stno
                 = 0\else .\number \stno \fi }\track{showlabel}{*}{#1}}
  \def \lcite #1{(#1\track{showcit}{$\bullet$}{#1})}

  \newcount \bibno \bibno =0
  \def \newbib #1{\global\advance\bibno by 1 \edef #1{\number\bibno}}
  \def \bibitem #1#2#3#4{\smallskip \item {[#1]} #2, ``#3'', #4.}
  \def \references {
    \begingroup 
    \bigskip \bigskip \goodbreak 
    \eightpoint 
    \centerline {\tensc References}
    \nobreak \medskip \frenchspacing }

  \def \axiom #1\par {\sysstate {Axiom}{\sl }{#1}}
  \def\imply{\ \Rightarrow\ }
  \def\mycorr{generalized correspondence}
  \def\Mycorr{Generalized correspondence}
  \def\Ker{{\rm Ker}}
  \def\hm{Hilbert-module}
  \def\hb{Hilbert-bimodule}
  \def  \lhm#1{left--#1--\hm}
  \def  \rhm#1{right--#1--\hm}
  \def\hbm#1#2{#1--#2--\hb}
  \def \bm#1#2{#1--#2--bimodule}

  \def\H{{\cal H}}
  \def\V{{\cal V}}
  \def\EV{E_\V}
  \def\EH{E_\H}
  \def\eV{e_{_\V}}
  \def\eH{e_{_\H}}

  \def\a{\alpha}
  \def\X{{\curly{X}}}   \def\Xs{{\hbox{\rssmall X}}}
  \def\TX{\T_\Xs}
  \def\OX{\O_\Xs}
  \def\XAT{\X\!\!_{\a,\Tr}}
  \def\XVH{\X\!_{\V,\H}}
  \def\Y{{\curly{Y}}}
  \def\J{{\cal J}}
  \def\I{{\cal I}}
  \def\S{{\cal S}}
  \def\T{{\cal T}}
  \def\TAX{\T(A,\X)}
  \def\O{{\cal O}}
  \def\M{{\cal M}}
  \def\D{{\cal D}}
  \def\L{\curly{L}}
  \def\R{\curly{R}}
  \def\Tr{{\cal L}}
  \def\Link{{\cal L}}
  \def\K{\curly{K}}
  \def\KV{\curly{K}\!_\V}
  \def\KH{\curly{K}\!\!_\H}
  \def\KL{\curly{K}\!\!_\ell}
  \def\KR{\curly{K}\!\!_r}
  \def\compos{\mathop{\scriptstyle \circ}} 
  
  \def\antiiso{\mathop{{\buildrel {op} \over \sim}}}
  \def\ip#1#2{\langle #1,#2 \rangle}
  \def\IP#1#2{\Big\langle #1,#2 \Big\rangle}
  \def\lp#1#2{\ip{#1}{#2}_\ell}
  \def\rp#1#2{\ip{#1}{#2}_r}
  \def\RP#1#2{\IP{#1}{#2}_r}
  \def\*{\otimes}
  \def\.{\odot}
  
  \def\[{\Big\|}
  \def\]{\Big\|}
  \def\half{^{1/2}}
  \def\({\left(}
  \def\){\right)}
  
  \def\tri(#1,#2,#3)(#4,#5,#6){[#1=#4][#2=#5][#3=#6]}%



  \newbib\Blackadar
  \newbib\BMS
  \newbib\Raeburn
  \newbib\Deaconu
  \newbib\TPA
  \newbib\endo
  \newbib\tower
  \newbib\vershik
  \newbib\Fowler
  \newbib\WataTwo
  \newbib\Kasparov
  \newbib\Katsura
  \newbib\Muhly
  \newbib\Pimsner
  \newbib\Rieffel
  \newbib\RieffelMorita
  \newbib\Takesaki
  \newbib\Watatani
  \newbib\Zettl

  \def\titletext{INTERACTIONS}

  \Headlines
  {\titletext}
  {R.~Exel}

  \null\vskip -1cm
  \centerline{\titlefont\titletext}

  \bigskip
  \centerline
  {\tensc R.~Exel\footnote{*}{\eightrm Partially supported by CNPq.}}

  \medskip\centerline
  {University of Santa Catarina}

  \medskip\centerline
  {September 2004}

  \bigskip \bigskip
  \midinsert\narrower

  \noindent{\tensc abstract.} Given a C*-algebra $B$, a closed
*-subalgebra $A\subseteq B$, and a partial isometry $S$ in $B$ which
\stress{interacts} with $A$ in the sense that
  $  \S^*a\S = \H(a)\S^*\S
  $ and $ 
  \S a\S^* = \V(a)\S\S^*,
  $
  where $\V$ and $\H$ are positive linear operators on $A$, we derive
a few properties which $\V$ and $\H$ are forced to satisfy.  Removing $B$
and $\S$ from the picture we define an \stress{interaction} as being a
pair of maps $(\V,\H)$ satisfying the derived properties.  Starting
with an abstract interaction $(\V,\H)$ over a C*-algebra $A$ we
construct a C*-algebra $B$ containing $A$ and a partial isometry $S$
whose \stress{interaction} with $A$ follows the above rules.
  We then discuss the possibility of constructing a \stress{covariance
algebra} from an interaction.  This turns out to require a
generalization of the notion of correspondences (also known as Pimsner
bimodules) which we call a \stress{\mycorr}.  Such an object should be
seen as an usual correspondence, except that the inner-products need
not lie in the coefficient algebra.  The covariance algebra is then
defined using a natural generalization of Pimsner's construction of the
celebrated Cuntz-Pimsner algebras.
  \endinsert

  \section{Introduction}
  In \cite{\endo} we have introduced a crossed-product construction
(see also \cite{\tower} and \cite{\vershik}) given an endomorphism
$\alpha$ of a C*-algebra $A$ and a transfer operator $\Tr$, attempting
to improve and extend previous constructions.  Since then a wealth of
interesting examples have been discovered, notably Watatani's study of
polynomial maps on subsets of the Riemann sphere \cite{\WataTwo} and
Deaconu's generalization of our construction to semigroups of
endomorphisms and his study of higher rank graphs \cite{\Deaconu} (see
also \cite{\Fowler}).

   Recall that, among other things, the crossed-product referred to
above involves a partial isometry $\S$ such that
  $$
  \S a = \alpha(a) \S
  \and
  \S^*a\S = {\cal L}(a)
  \for a\in A.
  $$

  If one looks at these relations from a purely aesthetical point of
view he or she will likely be bothered by their asymmetry.  For
instance, $\S$ is allowed to jump from left to right but not
vice-versa, that is, there is no relation telling us what to do with
an expression of the form $a\S$.

Perhaps the cause for this asymmetry is the fact that time
evolution is often irreversible.  In order to explain this recall that 
while $\a$ can be interpreted as accounting  for the future\fn{In Dynamical Systems
one often thinks of the given map as representing the time evolution
and hence its positive iterates represent the ``future''.} evolution
of the system under analysis, the transfer operator is needed to
account for the past.  In the case of \stress{irreversible systems},
namely systems in which $\a$ is not invertible, it is natural that
time evolution should behave quite differently depending on whether we
are moving forward or backward in time.

If we are speaking of a classical irreversible system, say a
continuous (non-invertible) map
  $$
  T : X \to X,
  $$
  where $X$ is a compact space, the trouble with accounting for
the past corresponds to the fact  that a  point $x$ in $X$ may have more
than one pre-image under $T$.  Yet the set $T\inv(x)$ of all
pre-images is always well defined and in some cases one can attach a
probability distribution to $T\inv(x)$ representing a guess as to what
the past looked like.  By looking at an already extinguished camp fire
it is impossible to tell how did it look like the night before but it
is sometimes possible to guess!

Suppose that one indeed is given a probability distribution $\mu_x$ on
$T\inv(x)$ for each $x$ in $X$. Then, given an observable, i.e.~a
continuous scalar valued  function $f$ on $X$, one may define
  $$
  \Tr(f)\calcat x = \int_{T\inv(x)} f(y)\,d\mu_x(y)
  \for x\in X.
  $$
  In this way $\Tr(f)$ represents the expected value of the observable
$f$ one unit of time into the past.
  Supposing that $\Tr(f)$ is continuous for every $f$ one checks without
difficulty that $\Tr$ is a transfer operator for the endomorphism
$\a$ of $C(X)$ defined by
  $$
  \a : f\in C(X) \mapsto f\compos\a \in C(X).
  $$

In this paper we wish to take a first step toward the study of systems
whose future behavior presents the same degree of uncertainty as its
past.

Inspired by our previous crossed-product construction we postulate
that our given algebra of observables $A$ should be embedded in a
larger algebra $B$ (roughly playing the role of the crossed-product)
containing a partial isometry $\S$ which governs time evolution.  Time
evolution itself will be thought of as the interaction between $A$ and
$\S$, meaning the commutation relations between $\S$ and the elements
of $A$.  These commutation relations will be required to have the form
  $$
  \S^*a\S = \H(a)\S^*\S
  \and
  \S a\S^* = \V(a)\S\S^*,
  \eqno{(\dagger)}
  $$
  for all $a\in A$, where $\V$ and $\H$ are positive linear maps
defined on $A$.  One should think of these maps as corresponding to
the past and future evolution.  Which is which is not an issue since
the situation will be absolutely symmetric.

Starting with as few assumptions as possible we will derive properties
of these maps, thus motivating the definition of an
\stress{interaction}: a pair of maps $(\V,\H)$ satisfying the derived
properties.  In particular we will show that 
  $$
  \H(a)\S^*\S = \S^*\S\,\H(a)
  \and
  \V(a)\S\S^* = \S\S^*\V(a),
  $$
  thus eliminating any  asymmetry in $(\dagger)$ which the reader might
have suspected.

  A rule to be followed here is to restrict one's attention to
statements which refer exclusively to $A$, as opposed to the larger
algebra $B$ or the partial isometry $\S$.  Once the definition is
given we will fix an interaction $(\V,\H)$ over a unital C*-algebra
$A$ and, ignoring any previous information relative to $B$ or $\S$, we
will reconstruct an instance of $B$ and $\S$ in such a way that
$(\dagger)$ holds.

The crucial step in this reconstruction is a certain ternary ring of
operators, TRO for short \cite{\Zettl} (see also \cite{\TPA}).  There
are several equivalent ways to define a TRO but perhaps the quickest
one is that a TRO is a closed subspace $\Y$ of a C*-algebra such that
$\Y\Y^*\Y\subseteq \Y$.  There is also an axiomatic definition in
which the ternary operation
  $$
  (x,y,z) \in \Y\times\Y\times\Y \mapsto xy^*z \in\Y
  $$
  plays a predominant role.  In the present situation in which $A$ is
assumed to be a subalgebra of $B$, where $\S$ also lies, we take
$\Y$ to be the closed linear span of $A\S A$.  Given $a,b,c,d,e,f\in
A$ observe that
  $$
  (a\S b)(c\S d)^*(e\S f) =
  a\V(bd^*)\ \S\ \H(c^*e) f,
  $$
  which shows that $\Y$ is in fact a TRO.
  
  If one is to reconstruct $\Y$ from the abstract information
contained in the given interaction $(\V,\H)$, the natural thing to do
is to let $\X$ be the tensor product (over the complexes) of $A$ by
itself, that is, 
  $$ 
  \X = A\*A,
  $$
  (later we will use ``$\.$'' for  the algebraic tensor
product over the complex numbers) and equip $\X$ with the ternary
operation 
  $$
  [\,\cdot\,,\,\cdot\,,\,\cdot\,]: \X\times\X\times\X \to \X
  $$
  defined by 
  $$
  [a\* b,c\* d,e\* f] =
  a\V(bd^*)\,\*\,\H(c^*e) f.
  $$
  Even if this is what we eventually accomplish we have found that it is
better to recover our TRO in the guise of a {\hb}.  Observe that every
{right--\hm} is a TRO relative to the ternary operation
  $$
  [\xi,\eta,\zeta] = \xi \<\eta,\zeta\>.
  $$
  In fact a TRO can be considered as a {\hb} in which the
  coefficient algebras are not explicit.
  Thus, supposing one can guess the coefficient algebras, describing a
TRO as a {\hb} is sometimes easier.  In the present case the right
guess is to take the reduced C*-basic constructions $\KV$ and $\KH$
\cite{\Watatani} relative to the conditional expectations
  $$
  \EV = \V\circ\H
  \and
  \EH = \H\circ\V,
  $$
  respectively.  Our major effort is thus directed at the task of
giving a \hbm{$\KV$}{$\KH$} structure to a certain completion of
$A\*A$. 
  The idea is to adopt the following ``dictionary''
  $$
  \matrix{\S & \to & \* \cr
          \S\S^* &  \to & \eV\cr
          \S^*\S &  \to & \eH}
  $$
  where $\eV$ and $\eH$ are the canonical idempotents in $\KV$ and
$\KH$, respectively.  Therefore the easily verifiable 
formulas 
  \begingroup
  \parindent = 30pt
  \izitem
  \zitem $(a\S b)(c\S d)^* = a\V(bd^*)\S\S^* c^*,$
  \zitem $(a\S b)^*(c\S d) = b^*\H(a^*c)\S^*\S d,$
  \zitem $(a\S\S^*b)(c\S d) = a \S \H(bc)d,$
  \zitem $(a\S b) (c\S^*\S d) = a\V(bc)\S d,$

  \medskip\noindent  translate to
  \izitem
  \Item{(i')} $\lp{a\* b}{c\* d\,} = a\V(bd^*)\eV c^*,$
  \Item{(ii')} $\rp{a\*b}{c\* d} = b^*\H(a^*c)\eH d,$
  \Item{(iii')} $(a\eV b)(c\*d) = a \* \H(bc)d,$
  \Item{(iv')} $(a\* b) (c\eH d) = a\V(bc)\* d,$
  \endgroup

  \medskip\noindent hence suggesting (i') the $\KV$--valued
inner-product, (ii') the $\KH$--valued inner-product, (iii') the
left--$\KV$--module structure, and (iv') the right--$\KH$--module
structure.
  
  Since the {\hm} we are after is in fact a \stress{completion} of
$A\*A$, in order to extend the above formulas beyond the algebraic
tensor product we must verify certain inequalities, the proof of which
turns out to a major headache.  It would therefore be interesting to
find shorter proofs to the existence of the above {\hm} structure.

Fortunately, after the TRO/{\hb} is constructed, it is easy to find an
algebra $B$ containing $A$ and a partial isometry $\S$ satisfying the
predicted commutation relations.

One is next faced with the question of extending the construction of
the crossed-product.  In the context of an endomorphism and a transfer
operator the crossed-product is the Cuntz-Pimsner algebra of a certain
Hilbert-bimodule (also known as a Pimsner bimodule \cite{\Pimsner} or
a correspondence \cite{\Muhly}, see also \cite{\Raeburn}).  

However a somewhat unexpected but very interesting novel feature
appears: the correspondence which seems to be the natural replacement
in the general case is not quite a correspondence!  Rather it is the
TRO $\X$ constructed above which is an $A$--$A$-bimodule but not a
\rhm{$A$}.  In other words, there is no $A$--valued inner-product
(although there is an inner-product taking values in a larger
algebra).  We call such an object a {\mycorr} (for the lack of a
better name).  Precisely speaking a {\mycorr} is a TRO which is at
the same time an $A$--$A$-bimodule such that
  \medskip\item{$\bullet$} $[\xi,a\eta,\zeta] =   [\xi,\eta,a^*\zeta]$, and
  \medskip\item{$\bullet$} $[\xi,\eta a,\zeta] =   [\xi
  a^*,\eta,\zeta]$,
  \medskip\noindent for all $a\in A$, and $\xi,\eta,\zeta\in\X$.

  As already mentioned a TRO is a {\hb} in which the coefficient
algebras are not explicit.  Nevertheless one may dig out the
coefficient algebras as follows: let  $\L(\X)$ and $\R(\X)$ be the
algebras formed by the left and right operators
\scite{\TPA}{4.3}, respectively, and let $\KL(\X)\subseteq\L(\X)$ and
$\KR(\X)\subseteq\R(\X)$ be the ideals of generalized compact
operators.  Then given any pair of C*-algebras $(L,R)$ with
  $$
  \KL(\X)\subseteq L\subseteq \L(\X)
  \and
  \KR(\X)\subseteq R\subseteq \R(\X),
  $$
  one sees that $\X$ is an \hbm{$L$}{$R$}.  The main difference
between the classical and generalized versions of correspondence is as
follows:  consider the *-anti-homomorphism
  $
  \rho: A \to \R(\X)
  $
  given my 
  $$
  \rho(a): \xi\in\X \mapsto \xi a \in\X.
  $$
  In case $\X$ is a classical correspondence one may easily prove that
every compact right operator coincides with $\rho(a)$ for some $a$ in
the ideal generated by $\rp\X\X$ and conversely, $\rho(a)$ is a
compact right operator for every $a$ in said ideal.  In particular
$\KR(\X) \subseteq \rho(A)$.  In the case often considered in which
$\X$ if full as a Hilbert-module (namely when $A$ is spanned by
$\rp\X\X$) one may prove that $\rho$ establishes an isomorphism
between $A$ and $\KR(\X)$ \scite{\BMS}{1.10}.

The major novel feature of a {\mycorr} is that the relationship
between the right action of $A$ on $\X$ and $\KR(\X)$ is not so
stringently restricted.

At this point our development becomes admittedly experimental.  Guided
by the fact that the right action on a generalized correspondence
displays properties similar to the left action on a usual
correspondence (namely a loose relationship to compact operators) we
attempt to treat both the left and right actions of $A$ on an equal
footing and we thus tentatively give a definition of
\stress{crossed-product} or \stress{covariance algebra}, denoted
$C^*(A,\X)$, using Pimsner's idea twice, namely \stress{moding out the
redundancies} on both sides.  See condition (4) of
\scite{\Pimsner}{Definition 3.8}.

This raises many more questions that we attempt to answer, in
particular we leave open the question as to whether the canonical
embeddings of $A$ and $\X$ in $C^*(A,\X)$ are injective (see below
for more questions).

We are nevertheless able to show that our generalized theory includes
the Cuntz-Pimsner algebra construction for the case of a usual
correspondence and the crossed-product construction in case of a pair
$(\a,\Tr)$ consisting of an endomorphism $\a$ and a transfer operator
$\Tr$.

  \section {Motivation}
  \label\MotivationSection
  In order to motivate our choice of axioms let $A$ be a closed
*-subalgebra of a C*-algebra $B$ and let $\S\in B$ be a partial
isometry.  

The main assumption we will adopt is related to commutation relations
between $\S$ and the elements of $A$.  Namely we will suppose that for
every $a\in A$ there exist $b,c\in A$ such that $\S a \S^* =
b\S\S^*$, and $\S^*a\S = c\S^*\S$.  We would also like that $b$ and
$c$ have a functional dependency on $a$ so we will actually assume:

  \axiom
  \label \CovarianceAxiom
  There are positive bounded linear maps
  $
  \V,\H:A\to A,
  $
  such that
  $$
  \S^*a\S = \H(a)\S^*\S
  \and
  \S a\S^* = \V(a)\S\S^*
  \for a\in A.
  $$

  The first immediate consequence of \lcite{\CovarianceAxiom} is the
following simple:

  \state Proposition 
  \label \Commute
  Every $x\in \H(A)$ commutes with $\S^*\S$ and
every $y\in\V(A)$ commutes with $\S\S^*$.

  \proof For every $a\in A$ we have
  $$
  \H(a)\S^*\S = \S^*a\S = (\S^*a^*\S)^* =
  \big(\H(a^*)\S^*\S\big)^* =
  \S^*\S\H(a^*)^* =
  \S^*\S\H(a),
  $$
  proving the first statement.  The second one is proven similarly.
  \proofend

  A trivial example of \lcite{\CovarianceAxiom} is obtained by taking any
pair of maps $(\V,\H)$ and choosing $\S=0$.  Since this is clearly
undesirable, we will avoid this situation by adopting the following
nondegeneracy property:

  \axiom 
  \label \NonDegeneracyAxiom
  \izitem
  \zitem
  For any $x$ and $y$ in $C^*(\V(A))$ (the C*-algebra generated\fn{We
will see later that our assumptions imply that $\V(A)$ is a
C*-algebra.} by $\V(A)$), one has that
  $$
  x \S = y \S \imply x = y.
  $$
  \zitem For any $x$ and $y$ in $C^*(\H(A))$,
  $$
  \S x = \S y \imply x = y.
  $$

  This is related to the uniqueness of $\V$.  In fact, if $\tilde\V$
is another map satisfying \lcite{\CovarianceAxiom}, then
  $$
  \tilde\V(a) \S =
  \tilde\V(a) \S\S^*\S =
  \S a\S^*\S =
  \V(a) \S\S^*\S =
  \V(a) \S.
  $$
  Thus, at least if $\tilde\V(a)$ lies in $C^*(\V(A))$, we must have
$\tilde\V(a) = \V(a)$ by \lcite{\NonDegeneracyAxiom.i}. Obviously a similar
reasoning applies to $\H$.

Given $a$ in $A$ observe  that
  $$
  \V(a)\S =
  \V(a)\S\S^*\S =
  \S\S^*\V(a)\S =
  \V(\H(\V(a)))\S\S^*\S =
  \V(\H(\V(a)))\S,
  $$
  and one may similarly prove that
  $$
  \S \H(a) =
  \S \H(\V(\H(a))).
  $$

  Combining this with \lcite{\NonDegeneracyAxiom} we conclude:

  \state Proposition
  \label \TripleProduct 
  The maps $\V$ and $\H$ satisfy
  \izitem  
  \zitem $\V\H\V=\V$,  and
  \zitem $\H\V\H=\H$.

\medskip
Given $a$ and $b$ in $A$ let us compute $\S a \S^* b \S$ in the following two ways:
  $$
  \S a \S^* b \S =
  \S a \H(b)\S^* \S =
  \V\big(a \H(b)\big)\S \S^* \S =
  \V\big(a \H(b)\big)\S,
  $$
  while
  $$
  \S a \S^* b \S =
  \V(a) \S\S^* b \S =
  \V(a) \S\H(b) \S^* \S =
  \V(a) \V\big(\H(b)\big) \S.
  $$
  Using \lcite{\NonDegeneracyAxiom.i} this implies that
  $$
  \V(a \H(b)) = \V(a) \V(\H(b)).
  $$
  Since both $\V$ and $\H$ are positive, and hence preserve the involution, this implies that 
  $$
  \V(\H(a)b) =
  \V(\H(a)) \V(b),
  $$
  as well.  We thus have:
  
  \state Proposition
  \label \Homomorphism
  Let  $x,y\in A$.
  \izitem
  \zitem If either $x$ or $y$ belong to $\H(A)$, than $\V(xy) =
\V(x)\V(y)$.
  \zitem If either $x$ or $y$ belong to $\V(A)$, than $\H(xy) =
\H(x)\H(y)$.
  
  \proof  Point (i) was proved above while (ii) follows from computing 
$\S^* a \S b \S^*$ in two different ways.
  \proofend

We are now in a position to show the fact already  hinted at that
$\V(A)$ and $\H(A)$ are subalgebras.
  
  \state Proposition
  \label \TheSubalgebras
  Let $\EV = \V\compos \H$ and $\EH = \H\compos\V$.  Then
  \izitem
  \zitem $\V(A)$ and $\H(A)$ are  closed *-subalgebras of $A$,
  \zitem $\EV$ is a conditional expectation onto $\V(A)$,
  \zitem $\EH$ is a conditional expectation onto $\H(A)$.

  \proof Given $a,b\in A$ we have that
  $$
  \V(a)\V(b) =
  \V(a)\V(\H(\V(b))) =
  \V(a\H(\V(b))),
  $$
  which shows that $\V(A)$ is a subalgebra of $A$.  Since $\V$ is
positive, and hence preserves the involution, we have that $\V(A)$ is
closed under the involution.  By \lcite{\TripleProduct.i} it is easy to see that for any
$a\in A$ one has that 
  $$
  a\in\V(A) \ \Longleftrightarrow \ \V(\H(a)) = a,
  $$
  from where we deduce that $\V(A)$ is norm-closed.  In a similar way
one shows that $\H(A)$ is a closed *-subalgebra of $A$.

Addressing (ii) recall that both $\V$ and $\H$ are positive maps and
hence so is their composition.  By \lcite{\TripleProduct.i} we have
that $\EV$ is an idempotent map whose range is precisely $\V(A)$.
Moreover, given $a\in A$ and $b\in\V(A)$, say $b=\V(c)$, where $c\in
A$, we have by \lcite{\Homomorphism} that
  $$
  \EV(ab)=
  \V\Big(\H\big(a\V(c)\big)\Big) =
  \V\Big(\H(a)\H\big(\V(c)\big)\Big) \$=  
  \V\big(\H(a)\big)\V\Big(\H\big(\V(c)\big)\Big) =  
  \EV(a) \V(c) =
  \EV(a) b.
  $$
  By taking adjoints one shows that 
  $
  \EV(ba)=
  b \EV(a)
  $ 
  as well, thus proving that $\EV$ is in fact a conditional
expectation onto $\V(A)$.  The proof of (iii) is similar.
  \proofend

  Our next result shows that $\V$, as well as $\H$, are obtained by the
composition of a conditional expectation and a *-isomorphism.

  \state Proposition 
  \label \InverseIsos
  Denote by $\V_1$ the restriction of\/ $\V$ to
$\H(A)$, considered as a mapping
  $$
  \V_1: \H(A) \to \V(A),
  $$
  and similarly let 
  $$
  \H_1: \V(A) \to \H(A)
  $$
  be given by $\H_1=\H|_{\V(A)}$.  Then both $\V_1$ and $\H_1$ are
*-isomorphisms, each being the inverse of the other, and we have
  $
  \V = \V_1 \EH
  $ and $
  \H = \H_1  \EV.
  $

  \proof By \lcite{\Homomorphism} it is clear $\V_1$ and $\H_1$ are
homomorphisms which preserve the involution as a consequence of
positivity.  Given $x\in\H(A)$, say $x=\H(a)$, where $a\in A$, we have
that
  $$
  \H_1(\V_1(x)) =
  \H(\V(x)) =
  \H(\V(\H(a))) =
  \H(a) = 
  x,
  $$
  so that $\H_1\V_1$ is the identity map on $\H(A)$ and a similar
reasoning shows that $\V_1\H_1$ is the identity on $\V(A)$.  So
$\V_1=\H_1\inv$ as asserted.  For $a\in A$ observe that by
\lcite{\TripleProduct.i} one has
  $$
  \V_1(\EH(a)) = 
  \V(\H(\V(a))) =
  \V(a),
  $$ 
  proving that   $\V = \V_1 \EH$, while the proof that $\H = \H_1\EV$ follows along similar lines.
  \proofend

  \state Proposition
  \label \HomoOnB
  The maps
  \izitem
  \zitem
  $
  a\in \H(A) \mapsto \V(a)\S \S^* \in B,
  $
  and
  \zitem
  $
  a\in \V(A) \mapsto \H(a)\S^*\S \in B
  $
  \smallskip\noindent
  are isometric *-homomorphisms.

  \proof
  For $a,b\in \H(A)$ we have by \lcite{\Commute} and \lcite{\Homomorphism} that
  $$
  \big(\V(a)\S \S^*\big)\big(\V(b)\S \S^*\big) =
  \V(a)\V(b)\S \S^* =
  \V(ab)\S \S^*.
  $$
  We leave it for the reader to fill in the other requirements to show
(i) be a *-homomorphism.  The proof that (ii) is also a *-homomorphism
follows similarly.

  It is well known that injective *-homomorphism are automatically
isometric, so, in order to conclude the proof, it is enough to show
our maps to be injective.  For this observe that if $a\in \H(A)$ and
$\V(a)\S \S^*=0$ then also $0 = \V(a)\S \S^*\S =\V(a)\S$.  So, by
\lcite{\NonDegeneracyAxiom.i}, one has that $\V(a)=0$, and hence $a=0$ by
\lcite{\InverseIsos}.  Similarly one shows that (ii) is injective.
  \proofend

As a consequence we have:

  \state Corollary  
  \label\IsometryWithS
  For every $a\in A$ one has that 
  $$
  \|\V(a)\S \S^*\| = \|\V(a)\| \and
  \|\H(a)\S^*\S \| = \|\H(a)\|.$$

  \proof
  We prove only (i).  Given $a\in A$ one has that
  $$
  \|\V(a)\S \S^*\| =
  \|\V(\H(\V(a)))\S \S^*\| \={(\HomoOnB)}
  \|\H(\V(a))\| \={(\InverseIsos)}
  \|\V(a)\|.
  \proofend
  $$

  \section{Interactions}
  Motivated by the above considerations we will now introduce the main
object of interest.

  \definition
  \label \DefInteraction
  If $A$ is a C*-algebra then a pair $(\V,\H)$ of maps
  $$
  \V,\H:A\to A
  $$
  will be called an \stress{interaction} if
  \izitem   
  \zitem $\V$ and $\H$ are positive bounded linear maps,
  \zitem $\V\H\V=\V$, 
  \zitem $\H\V\H=\H$,
  \zitem $\V(xy)=\V(x)\V(y)$, if either $x$ or $y$ belong to $\H(A)$,
  \zitem $\H(xy)=\H(x)\H(y)$, if either $x$ or $y$ belong to $\V(A)$.

  \sysstate{Remark}{\rm}{Theorems \lcite{\TheSubalgebras} and
\lcite{\InverseIsos}, although proven in the slightly different
context of section \lcite{\MotivationSection}, were based exclusively on
hypotheses which hold for any interaction and hence they hold in
general.  We will therefore freely use their conclusions in what
follows.}

  One important consequence is:

  \state Corollary
  \label \CompletelyPositive
  Let $(\V,\H)$ be an interaction over the C*-algebra $A$.  Then both
$\V$ and $\H$ are completely positive and completely contractive maps.

  \proof
  It is well known that *-homomorphisms and conditional expectations
\scite{\Takesaki}{IV.3.4} share the properties in the statement.  Since
$\V = \V_1 \EH$ and $\H = \H_1 \EV$ by \lcite{\InverseIsos}, the
conclusion follows.
  \proofend

As for examples we have:

  \state Proposition
  \label \TransferInteraction
  Let $A$ be unital a C*-algebra, $\alpha$ be an endomorphism of $A$,
and $\Tr$ be a transfer operator for $\a$ such that $\Tr(1)=1$.  Then
$(\a,\Tr)$ is an interaction over $A$.

  \proof
  Recall \cite{\endo} that $\Tr$ is a positive linear map satisfying
$\Tr(a\a(b)) = \Tr(a)b$, for all $a,b\in A$.  Taking $a=1$ we conclude
that $\Tr\circ\a = id$.  With this one easily shows the validity of
\lcite{\DefInteraction.ii-iii}.  That \lcite{\DefInteraction.iv} holds
is obvious and speaking of \lcite{\DefInteraction.v} observe that if
$y=\a(b)\in\a(A)$ then
  $$
  \Tr(xy) =   
  \Tr(x\a(b)) =   
  \Tr(x)b =
  \Tr(x)\Tr(\a(b)) =
  \Tr(x)\Tr(y).
  $$
  Given that $\Tr$ is positive, and hence self-adjoint, we see that
the above holds if $x\in\a(A)$, instead.
  \proofend

  The first step in going from the abstract notion of interaction to
the concrete situation of section \lcite{\MotivationSection} is taken
by considering the following concept:

  \definition
  \label\DefCovarRep
  A \stress{covariant representation} for the interaction $(\V,\H)$ over the
C*-alge\-bra $A$ is a triple $(\pi,\S,B)$, where $B$ is a C*-algebra,
$\pi$ is a *-homomorphism of $A$ into $B$, and $\S$ is a partial
isometry in $B$ such that
  \izitem 
  \zitem $\S\pi(a)\S^* = \pi(\V(A))\S\S^*$, and
  \zitem $\S^*\pi(a)\S = \pi(\H(A))\S^*\S$.

Without any further restrictions it is very easy to construct a
covariant representation: one can take any *-homomorphism  $\pi:A\to B$
whatsoever and choose $\S$ to be the zero operator.  In order to stay
away from this trivial case we will pay greater attention to 
covariant representations possessing the following nondegeneracy
property (see also \lcite{\NonDegeneracyAxiom}):

  \definition
  \label\DefineNonDeg
  A covariant representation $(\pi,\S,B)$ is said to be
\stress{nondegenerated} if:
  \izitem
  \zitem  For any $x$ and $y$ in $\V(A)$  one has that
  $$
  \pi(x)\S = \pi(y)\S \imply x = y.
  $$
  \zitem For any $x$ and $y$ in $\H(A)$,
  $$
  \S \pi(x) = \S \pi(y) \imply x = y.
  $$

In order to check nondegeneracy it is enough to worry about
\lcite{\DefineNonDeg.i} for, supposing the latter holds, let
$a\in\H(A)$ be such that $S\pi(a)=0$. Then
  $$
  0 = \S \pi(a)\S^*\S=
  \pi(\V(a))\S\S^*\S =
  \pi(\V(a))\S.
  $$
  So $\V(a)=0$ and hence $a=\EH(a)=\H(\V(a))=0$.  We leave it for the
reader to show that \lcite{\DefineNonDeg.ii} implies
\lcite{\DefineNonDeg.i} as well.

Observe also that, if either $\EV$ or $\EH$ are faithful, as
conditional expectations, any covariant representation $(\pi,\S,B)$ must
be such that $\pi$ is one-to-one.  In fact, suppose that $\pi(a)=0$.
Then
  $$
  0 = \S\pi(a^*a)\S^*\S = 
  \pi(\V(a^*a))\S\S^*\S =
  \pi(\V(a^*a))\S,
  $$ 
  which implies that $\V(a^*a)=0$, and hence also
$\EH(a^*a)=\H(\V(a^*a))=0$.  Under the assumption that $\EH$ is
faithful one would then conclude that $a=0$.

Let us now briefly consider interactions over matrix algebras.  For
this suppose that $(\V,\H)$ is an interaction over a C*-algebra $A$
and let $n$ be a positive integer.  
Consider the maps
  $$
  \V_n,\H_n:M_n(A)\to M_n(A)
  $$
  obtained by applying $\V$ and $\H$, respectively, on all matrix entries.

It is elementary to verify that $(\V_n,\H_n)$ satisfies
\lcite{\DefInteraction.ii-v}.  By \lcite{\CompletelyPositive} we have
that $\V_n$ and $\H_n$ are positive and contractive and hence
\lcite{\DefInteraction.i} holds as well.
This proves:

  \state Proposition
  \label \MatrixInteraction
  Let $(\V,\H)$ be an interaction over a C*-algebra $A$.  Then
$(\V_n,\H_n)$ is an interaction over $M_n(A)$ for each positive
integer $n$.

The elementary aspects of interactions and covariant representations
discussed, we now want to worry about the existence of nondegenerated
covariant representations for a given interaction.  In order to obtain
inspiration for this task we will once again return to the concrete
situation of section \lcite{\MotivationSection}.

  \section{A ternary ring of operators}
  \label \ConcreteTROSection
  We now return to the hypotheses of our motivation section
\lcite{section \MotivationSection} and hence we assume that $A$ is a
closed *-subalgebra of a C*-algebra $B$ and that $\S$ is a partial
isometry in $B$ such that axioms \lcite{\CovarianceAxiom} and
\lcite{\NonDegeneracyAxiom} hold with respect to a given interaction
$(\V,\H)$.

  Let $\Y_0$ be the linear subspace of $B$ given by
  $$
  \Y_0 =
  \Big\{\sum_{i=1}^n a_i\S b_i : n\in\N,\ a_i,b_i\in A\Big\}.
  $$

  \state Proposition
  \label \ConcreteTRO
  Given
  $
  a_1,a_2,b_1,b_2,c_1,c_2\in A
  $
  we have that
  $$
  (a_1\S a_2)(b_1\S b_2)^*(c_1\S c_2) =
  a_1 \V(a_2 b_2^*)\S\H(b_1^*c_1)c_2.
  $$

  \proof  We have
  $$
  (a_1\S a_2)(b_1\S b_2)^*(c_1\S c_2) =
  (a_1\S a_2)(b_2^*\S^* b_1^*)(c_1\S c_2) \$=
  a_1 \V(a_2 b_2^*)\S\S^*\S\H(b_1^*c_1)c_2 =
  a_1 \V(a_2 b_2^*)\S\H(b_1^*c_1)c_2.
  \proofend
  $$

  This shows that $\Y_0$ is closed under the ternary operation
  $$
  (x,y,z) \mapsto xy^*z,
  $$
  and hence so is its closure, which we denote by $\Y$.
  In other words $\Y$ is a TRO (ternary ring of operators) in the sense
  of Zettl \cite{\Zettl}.
  In particular, both
  $$
  \KL := \Y\Y^*
  \and
  \KR := \Y^*\Y
  $$
  (by $\Y\Y^*$ we mean the \stress{closed linear span} of the set of
products $mn^*$, where $m,n\in\Y$, and similarly for $\Y^*\Y$)  are 
closed *-subalgebras of $B$, and $\Y$ is a \hbm{$\KL$}{$\KR$}
with the  left and right inner-products given by
  $$ 
  \lp mn = mn^* 
  $$
  and
  $$
  \rp mn = m^*n,
  $$
  for any $m,n\in\Y$.  As a notational device, let us write $\Y_r$
when we view $\Y$ as a \rhm{$\KR$} and by $\Y_\ell$
when $\Y$ is viewed as a \lhm{$\KL$}.

  It is also clear that $\Y$ is an \bm{$A$}{$A$}, however not
necessarily an \hbm{$A$}{$A$} as there is no reason for $\rp mn$ or
$\lp mn$ to lie in $A$.

  Nevertheless
  for every $a\in A$ and $m,n\in\Y_r$ one has that
  $$
  \rp{am}n =  
  m^*a^*n =
  \rp m{a^*n},
  $$
  so the map 
  $$
  \lambda(a):m\in\Y_r\mapsto am\in\Y_r
  $$
  is adjointable with adjoint given by $\lambda(a^*)$.  This says that
the correspondence
  $
  a \mapsto \lambda(a)
  $
  is a *-homomorphism from $A$ into the C*-algebra $\L(\Y_r)$ of all
adjointable operators on $\Y_r$.  Similarly, for every $a\in A$, the assignment
  $$
  \rho(a) : m\in \Y_\ell \mapsto ma \in \Y_\ell
  $$
  defines an element $\rho(a)\in\L(\Y_\ell)$ and  the correspondence
  $
  a \mapsto \rho(a)
  $
  is a *-anti-homomor\-phism from $A$ into $\L(\Y_\ell)$.

  We now wish to obtain a formula for the norm of an element $x$ of $\Y_0$
in terms of the maps $\V$ and $\H$.  For this write 
  $$
  x = 
  \sum_{i=1}^n a_i^*\S b_i
  $$
  (clearly there is no harm in using $a_i^*\S b_i$ where one might
expect $a_i\S b_i$, the reason for which will soon become clear), where
$n\in\N$ and $a_i,b_i\in A$.  Since $\|x\|^2 = \|x^*x\|$ we compute
  $$
  x^*x =
  \sum_{i,j=1}^n b_i^*\S^*a_i a_j^*\S b_j =
  \sum_{i,j=1}^n b_i^*\H(a_i a_j^*)\S^*\S b_j.
  $$
  If we view 
  $$
  a= \pmatrix{a_1 \cr a_2 \cr \vdots \cr a_n}
  $$
  as an $n{\times1}$ column matrix then $aa^*$ is a positive $n\times
n$ matrix.  By \lcite{\CompletelyPositive} it follows that
$\H_n(aa^*)$ is a positive element in the subalgebra
$M_n(\H(A))\subseteq M_n(A)$.  So there exists a matrix $c =
\big(c_{ij}\big)_{i,j=1}^n$ in that subalgebra such that
$\H_n(aa^*)=c^*c$.

Henceforth we will make extensive use of the maps $\V_n$ and $\H_n$ of
\lcite{\MatrixInteraction}.  However we will prefer to drop the
subscript ``$n$'' since the context will suffice to distinguish $\V$ and $\H$
from their matricial counterparts.

For all $i,j=1,\ldots,n$ we then have
  $
  \H(a_ia_j^*) =
  \sum_{k=1}^n c_{ki}^*c_{kj}.
  $
  So
  $$
  x^*x = 
  \sum_{i,j,k=1}^n b_i^*c_{ki}^*c_{kj}\S^*\S b_j \={(\Commute)}
  \sum_{i,j,k=1}^n b_i^*c_{ki}^*\S^*\S c_{kj}b_j \$=
  \sum_{k=1}^n \Big(\sum_{i=1}^n b_i^*c_{ki}^*\S^*\Big) 
               \Big(\sum_{j=1}^n \S c_{kj}b_j\Big) =
  (\S cb)^*(\S cb),
  $$
  where the correct interpretation of $\S cb$ above requires identifying $\S$
with the diagonal $n\times n$ matrix 
  $$
  \pmatrix{\S & 0 & \cdots & 0 \cr
                 0 & \S & \cdots & 0 \cr
                 \vdots & \vdots & \ddots & \vdots \cr
                 0 & 0 & \cdots & \S \cr} \in M_n(B).
  $$
  It follows that
  $\|x\|^2 =  \|(\S cb)^*(\S cb)\| = \|(\S cb)(\S cb)^*\|$
  which motivates our interest in the computation of
  $$
  (\S cb)(\S cb)^* =
  \S cbb^*c^*\S^* =
  \V(cbb^*c^*)\S\S^*.
  $$
  By the matricial version of \lcite{\IsometryWithS} we conclude that
  $$
  \|x\|^2 = 
  \|\V(cbb^*c^*)\| = 
  \|\V(c)\V(bb^*)\V(c^*)\| = 
  \|\V(c)\V(bb^*)\half \V(bb^*)\half\V(c^*)\| \$= 
  \|\V(bb^*)\half\V(c^*) \V(c)\V(bb^*)\half \| = 
  \|\V(bb^*)\half\V(c^*c)\V(bb^*)\half \| \$= 
  \|\V(bb^*)\half\V\big(\H(aa^*)\big)\V(bb^*)\half \| = 
  \|\V(bb^*)\half\V\big(\H(aa^*)\big)\half\V\big(\H(aa^*)\big)\half\V(bb^*)\half \| \$= 
  \|\V\big(\H(aa^*)\big)\half\V(bb^*)\half \|^2.
  $$

This gives the following:

  \state Proposition
  \label \NormComputation
  If $a_i,b_i\in A$, for $i=1,\ldots,n$ then
  $$
  \matrix{
  \big\|
  \sum_{i=1}^n a_i^*\S b_i
  \big\| & = & \big\| \V\big(\H(aa^*)\big)\half  \hfill
                      \V(bb^*) \half \big\| \cr\cr
         & = & \big\|\H(aa^*)\half\H\big(\V(bb^*)\big)\half \big\|
  }
  $$

  \proof The first identity was proved just above.  As for the second one
we claim that 
  $$
  \big\|\H(aa^*)\half\H\big(\V(bb^*)\big)\half \big\| \leq 
  \big\|\V\big(\H(aa^*)\big)\half \V(bb^*)\half \big\|.
  $$ 
  In fact,
  $$
  \big\|\H(aa^*)\half\H\big(\V(bb^*)\big)\half \big\| = 
  \big\|\H\Big(\V\big(\H(aa^*)\half\big)\Big)\H\big(\V(bb^*)\half\big) \big\| \$= 
  \big\|\H\Big(\V\big(\H(aa^*)\half\big)\ \V(bb^*)\half\Big) \big\| \leq
  \big\|\V\big(\H(aa^*)\half\big) \V(bb^*)\half \big\| \$=
  \big\|\V\big(\H(aa^*)\big)\half \V(bb^*)\half \big\|.
  $$ 
  This proves our claim and since the reverse inequality can be
similarly proved we conclude that equality holds.
  \proofend

  \section{From interactions to ternary ring of operators}
  \label\ReconstructionSection
  The reader might have noticed that, except for the introduction, our
sections have alternated between the concrete setup of section
\lcite{\MotivationSection} and the setup of abstract interactions.
Keeping up with this pattern we now fix an arbitrary interaction
$(\V,\H)$ over a C*-algebra $A$.

The major goal of this work is to construct a C*-algebra $B$
containing $A$ and a partial isometry $\S\in B$ satisfying
\lcite{\CovarianceAxiom}.
  In this section we take the significant intermediate step of
constructing the ternary ring of operators corresponding to the $\Y$
of the previous section.

  A ternary ring of operators is necessarily a {\hb} over
the algebra of generalized compact left and right operators
\scite{\TPA}{Section 4} and hence we may construct it as such.  This
is especially convenient since we already have a good hint as to what
the coefficient algebras should be: briefly returning to the concrete
situation of section \lcite{\MotivationSection} observe that for all $a\in
A$ we have
  $$
  \S\S^* a\, \S\S^* = \EV(a)\S\S^*
  \and
  \S^*\S\, a\, \S^*\S = \EH(a)\S^*\S.
  $$
  Doesn't this look suspiciously like the algebraic relations underlying the
basic construction of V.~Jones (see, e.g.~\cite{\Watatani}) for the
conditional expectations $\EV$ and $\EH$?  We will therefore denote by
$\KV$ and $\KH$ the reduced C*-basic constructions associated to these
conditional expectations, respectively, and we will try to reconstruct
our ternary ring of operators in the guise of a \hbm{$\KV$}{$\KH$}.

  The canonical idempotent in $\KV$ (resp.~$\KH$) will be
denoted by $\eV$ (resp.~$\eH$).  Given the relations above one should
think of $\eV$ as $\S\S^*$ and $\eH$ as $\S^*\S$.

In our next definition we use the symbol ``$\.$'' to means the
algebraic tensor product over the complex numbers.

  \definition
  \label \DefineInnerProds
  Let
  $
  \X_0 = A \. A
  $
  and let
  $$
  \lp{\,\cdot\,}{\,\cdot\,} : \X_0 \times \X_0 \to \KV
  \and
  \rp{\,\cdot\,}{\,\cdot\,} : \X_0 \times \X_0 \to \KH
  $$
  be the sesqui-linear functions defined  by
  $$
  \lp{a\.b}{c\.d\,} = a\V(bd^*)\eV c^*
  \and
  \rp{a\.b}{c\.d\,} = b^*\H(a^*c)\eH d.
  $$

It should be understood that   $\lp{\,\cdot\,}{\,\cdot\,}$ is linear
in the first variable and conjugate-linear in the second   one, while
the reverse applies for    $\rp{\,\cdot\,}{\,\cdot\,}$.

  We urge the reader to compute $(a\S b) (c\S d)^*$ and $(a\S b)^*
(c\S d)$ in the context of section \lcite{\MotivationSection} for motivation
of the above definition.

We do not claim that $\X_0$ is a {\hb} yet but a certain
completion of it will be given such a structure shortly.  We are
nevertheless already in a position to show the crucial positivity
axiom.

  \state Proposition
  \label \Positivity
  If $x\in \X_0$ then $\lp xx$ and $\rp xx$ are non-negative elements
of $\KL$ and $\KR$, respectively.

  \proof
  Write  $x = \sum_{i=1}^n a_i^*\. b_i$ so that
  $$
  \rp xx =
  \sum_{i,j=1}^n \rp {a_i^*\. b_i}{a_j^*\. b_j} =
  \sum_{i,j=1}^n b_i^*\H(a_ia_j^*)\eH b_j.
  $$
  Since $\H$ is completely positive we have that
$\big(\H(a_ia_j^*)\big)_{i,j} = \H(aa^*)$ is a non-negative element of
$M_n(\H(A))$.  It follows that that there exists $c\in M_n(\H(A))$
such that $\H(aa^*)=c^*c$, i.e., that
  $
  \H(a_ia_j^*) =
  \sum_{k=1}^n c_{ki}^*c_{kj},
  $
  for all $i$ and $j$.
  Therefore, recalling that $\eH$ commutes with the elements in the range
of $\EH$, i.e.~$\H(A)$, we have
  $$
  \rp xx =
  \sum_{i,j,k=1}^n b_i^* c_{ki}^*c_{kj} \eH b_j  =
  \sum_{i,j,k=1}^n b_i^* c_{ki}^* \eH  c_{kj} b_j  \$=
  \sum_{k=1}^n
    \Big(\sum_{i=1}^n  c_{ki} b_i\Big)^*\eH \Big(\sum_{j=1}^n c_{kj} b_j\Big) \geq0.
  $$
  A similar argument shows that $\lp xx \geq0$.
  \proofend

It is now time to prove the {\hm} version of the
Cauchy-Schwartz inequality.  Once again, although we do not yet have a
{\hm} we have enough to make the standard proof work.

  \state Lemma
  \label \CauchySchwartz
  Given $x,y\in \X_0$ one has
  \izitem 
  \zitem $\rp xy \rp yx \leq \rp xx \big\|\rp yy\big\|$, and
  \zitem $\lp xy \lp yx \leq \lp xx \big\|\lp yy\big\|$.
   
  \proof 
  Let
  $x = \sum_{i=1}^n a_i\. b_i$,
  $y = \sum_{i=1}^n c_i\. d_i$, and set
  $$
  h = \rp yx = \sum_{k,l=1}^n d_k^*\H(c_k^*a_l)\eH b_l.
  $$
  We wish to make sense of the product $y h$, even though $\X_0$ lacks
the structure of a right-module over $\KH$ and $h\in\KH$.  The idea is
to define an ad-hoc version of $y h$ by putting
  $$
  y' = 
  \sum_{i,k,l=1}^n c_i \V(d_i  d_k^*) \.  \H(c_k^*a_l)b_l.
  $$
  We claim that  for any $z\in A\.A$ one has that
  $$
  \rp{z}{y'} = \rp{z}{y}h.
  $$
  In order to verify our claim it is obviously enough to consider only
  the case $z=a\.b$, where $a,b\in A$.
  Given that, the left-hand side equals
  $$
  \rp{z}{y'} =   \rp{a\.b}{y'} = 
  \sum_{i,k,l=1}^n \rp{a\.b}{c_i \V(d_i  d_k^*) \.  \H(c_k^*a_l)b_l} \$=
  \sum_{i,k,l=1}^n b^* \H\big(a^* c_i \V(d_i  d_k^*)\big) \eH  \H(c_k^*a_l)b_l =
  \sum_{i,k,l=1}^n b^* \H(a^* c_i) \H\big(\V(d_i  d_k^*)\big) \eH  \H(c_k^*a_l)b_l,
  $$
  while the right-hand side equals
  $$
  \rp{z}{y}h = \rp{a\.b}{y}h =
  \sum_{i,k,l=1}^n \rp{a\.b}{c_i\.d_i} d_k^*\H(c_k^*a_l)\eH b_l \$=
  \sum_{i,k,l=1}^n b^*\H(a^* c_i)\eH d_i d_k^*\H(c_k^*a_l)\eH b_l =
  \sum_{i,k,l=1}^n b^*\H(a^* c_i)\EH\big( d_i d_k^*\H(c_k^*a_l)\big)\eH b_l \$=
  \sum_{i,k,l=1}^n b^*\H(a^* c_i)\H\big(\V( d_i d_k^*)\big)\H(c_k^*a_l)\eH b_l,
  $$
  proving the claim.  By taking adjoints it is clear that
  $\rp{y'}{z} = h^*\rp{y}{z}$, as well.

  By \lcite{\Positivity} we have for every real number $\lambda$ that
  $$
  0 \leq
  \rp{x - \lambda y'}{x - \lambda y'} =
  \rp xx - \rp x{\lambda y'} - \rp {\lambda y'}x + 
    \rp {\lambda y'}{\lambda y'} \$=
  \rp xx - \lambda \rp x{y}h - \lambda h^*\rp {y}x + 
    \lambda^2 h^*\rp{y}{y}h \$=
  \rp xx - 2\lambda h^*h  + \lambda^2 h^*\rp{y}{y}h.
  $$
  It is now well known that this implies that
  $$
  h^*h \leq \rp xx \big\|\rp yy\big\|,
  $$
  proving (i).  We leave it for the reader to adapt the above proof in
order to check (ii).
  \proofend

  It is now routine to prove that the expressions
  $$
  \|x\|_r:= \|\rp xx\|\half \and   \|x\|_\ell:= \|\lp xx\|\half 
  $$
  define seminorms on $\X_0$.
  The next result is inspired by \lcite{\NormComputation}.

  \state Proposition
  \label \NormComputationAbstract
  Let $x = \sum_{i=1}^n a_i^*\. b_i \in \X_0$, where $a_i,b_i\in A$.  
Viewing $a$ and $b$ as $n\times 1$  column matrices one has that
  $$
  \|x\|_r  =
  \|\H(aa^*)\half\H(\V(bb^*))\half\| =
  \|\V(\H(aa^*))\half\V(bb^*)\half\| =
  \|x\|_l.
  $$

  \proof 
  As in the proof of \lcite{\Positivity} let $c\in M_n(\H(A))$ be such
$\H(aa^*)=c^*c$.  We will actually suppose that $c= \H(aa^*)\half$.
Picking up from the last computation in the proof of
\lcite{\Positivity} we have that
  $$
  \rp xx =
  \sum_{k=1}^n\Big(\sum_{i=1}^n c_{ki} b_i\Big)^*\eH \Big(\sum_{j=1}^n
    c_{kj} b_j\Big) =
  (cb)^*\eH cb,
  $$
  where, in the last term above, we identify $\eH$ with  the diagonal $n\times n$ matrix 
  $$
  \pmatrix{\eH & 0 & \cdots & 0 \cr
                 0 & \eH & \cdots & 0 \cr
                 \vdots & \vdots & \ddots & \vdots \cr
                 0 & 0 & \cdots & \eH \cr} \in M_n(\KH).
  $$
  It follows that
  $$
  \| \rp xx \| =
  \|b^*c^*\eH cb\| =
  \|\eH cbb^*c^*\eH\| =
  \|\EH(cbb^*c^*)\eH\| = \ldots 
  $$
  Recall from \scite{\Watatani}{2.1.4(2)} that for an element
$x\in\H(A)$ one has that $\|x\eH\|=\|x\|$.  This is so in spite of the
fact that our conditional expectations are not necessarily faithful,
while the conditional expectations in \cite{\Watatani} are required to
possess this property.  Fortunately the proof of the result referred
to does not depend on faithfulness.

So the above equals
  $$
  \ldots = 
  \|\EH(cbb^*c^*)\| =
  \|c\EH(bb^*)c^*\| =
  \|c\EH(bb^*)\half\|^2 \$=
  \|\H(aa^*)\half\H(\V(bb^*))\half\|^2.
  $$
  A similar reasoning shows that 
  $$
  \| \lp xx \| =
  \|\V(\H(aa^*))\half\V(bb^*)\half\|^2,
  $$
  and the conclusion follows as in \lcite{\NormComputation}.
  \proofend

  From now on we will consider $\X_0$ as a semi-normed space relative
to either $\|\cdot\|_\ell$ or to $\|\cdot\|_r$, which we henceforth
denote simply by $\|\cdot\|$.  Moding out by the vectors of length
zero and completing we get a Banach space which we will denote by
$\X$.
  If $a,b\in A$ then the canonical image of the element $a\.b$ in $\X$
will be denoted by $a\*b$.
By \lcite{\CauchySchwartz} we have that
  $$
  \|\lp xy \| \leq \|x\|\, \|y\| \and 
  \|\rp xy \| \leq \|x\|\, \|y\|,
  $$
  for every $x,y\in\X_0$ and hence we may continuously extend the two
inner-products to sesqui-linear maps
  $$
  \lp{\,\cdot\,}{\,\cdot\,} : \X \times \X \to \KV
  \and
  \rp{\,\cdot\,}{\,\cdot\,} : \X \times \X \to \KH.
  $$
  For all $a,b,c,d\in A$ one therefore has the ``$\*$'' versions of
\lcite{\DefineInnerProds}:
  $$
  \lp{a\*b}{c\*d\,} = a\V(bd^*)\eV c^*
  \and
  \rp{a\*b}{c\*d\,} = b^*\H(a^*c)\eH d.
  \eqno{(\seqnumbering)}  
  \label \InnerProduct
  $$ 

  \state Proposition
  \label \Balanced
  Let $a,b,c\in A$.  
  \izitem
  \zitem  If $c\in\V(A)$ then $ac\*b=a\*\H(c)b$.
  \zitem  If $c\in\H(A)$ then $a\*cb=a\V(c)\*b$.

  \proof
  In order to prove  (i) we need to verify that
  $\|ac\*b- a\*\H(c)b\| = 0$, for which it suffices to check that
  $$
  \rp{x\*y}{ac\*b - a\*\H(c)b} = 0,
  $$
  for every $x,y\in A$.
  We have
  $$
  \rp{x\*y}{ac\*b} =
  y^*\H(x^*ac)\eH b =
  y^*\H(x^*a)\H(c)\eH b \$=
  y^*\H(x^*a)\eH \H(c)b =
  \rp{x\*y}{a\*\H(c)b}.
  $$
  The proof of (ii) follows as above with the computation
  $$
  \lp{a\*cb}{x\*y} =
  a\V(cby^*)\eV x^* =
  a\V(c)\V(by^*)\eV x^* =
  \lp{a\V(c)\*b}{x\*y}.
  \proofend
  $$

We now wish to give the first step toward the proof of the existence
of a \hbm{$\KV$}{$\KH$} structure on $\X$.  Since the computations
involved are rather long we will concentrate on the
\rhm{$\KH$} structure for the time being.  For motivation one
should observe that, in the context of section
\lcite{\MotivationSection}, one has that
  $$
  (x\S y)(a\S^*\S b) = x\V(ya)\S b,
  \eqno{(\seqnumbering)}   
  \label \SuggestRModStru
  $$
  for all $a,b,x,y\in A$.

\definition Given $a,b\in A$ we will denote by $R_{a,b}$ the operator on $\X_0$
given by
  $$
  R_{a,b}(x\.y) = x\V(ya)\.b,
  $$
  for every $x,y\in A$.

The next technical lemma is intended to provide the key inequality
which will enable us to extend this to a right--$\KH$--module structure
on $\X$.  The length of its proof perhaps calls for a shorter argument
which we unfortunately have not yet found.

  \state \label\TechLemma
  Lemma Let  $N\in\N$ and for  each $n=1,\ldots,N$, let
$a_n,b_n\in A$.  Then
  $$
  \Big\|\RP \xi{\sum_{k=1}^N R_{a_k^*,b_k} (\eta)}\Big\| \leq \|\xi\|\, \|\eta\|\, \|\phi\|
  \for \xi,\eta\in \X_0,
  $$
  where $\phi$ is the element of $\KH$ defined by $\phi  = \sum_{k=1}^N a_k^* \eH b_k$.

  \proof Write  $$
  \xi  = \sum_{i=1}^N x_i^* \. y_i
  \and
  \eta = \sum_{j=1}^N z_j^* \. w_j, 
  $$
  where $x_n,y_n,z_n,w_n\in A$ for every $n=1,\ldots,N$.  In fact the
number of summands in the definition of $\xi$ and $\eta$ need not
necessarily be $N$, but we may clearly adjust to this situation by padding
the above sums with zeros in order to match the number of summands.

  \newcount\eqcntr \eqcntr=0
  \def\eqmark#1{\global \advance\eqcntr by 1 
    \edef\a{\number \secno .\number\stno.\number\eqcntr}
    \eqno {(\a)}
    \global \edef #1{\a}}

We will initially deal with the simplified case in which $N=1$ and, in
order to lighten  notation, we will write $x,y,z,w,a,b$ in place of
$x_1,y_1,z_1,w_1,a_1,b_1$.
  We therefore have that
  $$
  \rp \xi {R_{a^*,b}(\eta)} =
  \rp {x^*\.y}{z^*\V(wa^*)\.b} =
  y^*\H\big(x z^*\V(wa^*)\big)\eH b \$=
  y^*\H(x z^*)\eH \H\big(\V(wa^*)\big)b =
  y^*\H(x z^*)\eH \EH(wa^*)b =
  u^*\eH v,
  \eqmark\Contas
  $$
  where $u$ and $v$ are defined by
  $u  =  \H(zx^*)y$, and $v = \EH(wa^*)b$.
  We then conclude from \scite{\Watatani}{2.2.4} (again this does not
depend on faithfulness) that
  $$
  \|\rp \xi{R_{a^*,b}(\eta)}\|^2 =  
  \|u^*\eH v\|^2=
  \|\EH(uu^*)\half\EH(vv^*)\half\|^2 \$=
  \|\EH(vv^*)\half\EH(uu^*)\EH(vv^*)\half\|.
  \eqmark\EqOne
  $$
  Observe that
  $$
  \EH(uu^*) =
  \EH\big( \H(zx^*)y y^*\H(x z^*)\big) =
  \H(zx^*)\EH(y y^*)\H(x z^*) \$=
  \H\V\H(zx^*)\ \EH(y y^*)\ \H\V\H(x z^*) =
  \H\EV(zx^*)\ \H\V(y y^*)\ \H\EV(x z^*) \$=
  \H\Big(\EV(zx^*)\ \V(y y^*)\ \EV(x z^*)\Big) =
  \H\Big( \EV\big(zx^*\V(y y^*)\half\big) \
    \EV\big(\V(y y^*)\half x z^*\big) \Big).
  $$
  Recall from \scite{\Rieffel}{2.9} that every conditional expectation $E$
satisfies the inequality 
  $$
  E(cd^*)E(dc^*)\leq E(cc^*)\ \|E(dd^*)\|.
  \eqmark\EqTwo
  $$
  If we apply this with $c=z$, and $d=\V(y y^*)\half x$, we conclude that
  $$
  \EH(uu^*) \leq
  \H\Big(
    \EV(zz^*)\ 
    \big\| \EV\big( \V(y y^*)\half x x^* \V(y y^*)\half
    \big)\big\|\Big).
    \eqmark\EqThree
  $$
  The term between the double bars above satisfies
  $$
  \EV\big( \V(y y^*)\half x x^* \V(y y^*)\half \big) =
  \V(y y^*)\half \EV(x x^*) \V(y y^*)\half,
  $$
  and its norm is given by
  $$
  \big\| \V(y y^*)\half \EV(x x^*) \V(y y^*)\half \big\|=
  \big\| \V\big(\H(x x^*)\big)\half \V(y y^*)\half \big\|^2 \={(\NormComputationAbstract)}
  \|x^*\.y\|^2.
  $$
  By \lcite{\EqThree} we have 
  $$
  \EH(uu^*) \leq
  \|x^*\.y\|^2\ \H\big(\EV(zz^*)\big) =
  \|x^*\.y\|^2\ \H(zz^*).
  $$
  Returning to \lcite{\EqOne} we then have 
  $$
  \|\rp \xi{R_{a^*,b}(\eta)}\|^2 =  
  \|\EH(vv^*)\half\EH(uu^*)\EH(vv^*)\half\| \$\leq
  \|x^*\.y\|^2\ \|\EH(vv^*)\half\H(zz^*)\EH(vv^*)\half\| \$=
  \|x^*\.y\|^2\ \|\H(zz^*)\half\EH(vv^*)\H(zz^*)\half\|.
  \eqmark\EqFour
  $$
  Dealing with the term $\EH(vv^*)$ occurring above we have that
  $$
  \EH(vv^*) =
  \EH\big(\EH(wa^*)b b^* \EH(aw^*)\big) =
  \EH(wa^*) \EH(b b^*) \EH(aw^*) \$=
  \EH\big(wa^*\EH(b b^*)\half\big) \EH\big(\EH(b b^*)\half aw^*\big).
  $$
  Applying \lcite{\EqTwo} again with $c = w$, and $d = \EH(b b^*)\half a$, gives
  $$
  \EH(vv^*) \leq
  \EH(ww^*)\ \big\|\EH\big( \EH(b b^*)\half a a^* \EH(b b^*)\half\big)\big\| \$=
  \EH(ww^*)\ \big\|\EH(b b^*)\half \EH(a a^*) \EH(b b^*)\half\big\| \$=
  \EH(ww^*)\ \big\|\EH(a a^*)\half \EH(b b^*)\half\big\|^2 =
  \EH(ww^*)\ \|a^*\eH b\|^2.
  $$
  With \lcite{\EqFour} we now conclude that
  $$
  \|\rp \xi{R_{a^*,b}(\eta)}\|^2 \leq
  \|x^*\.y\|^2\ \|\H(zz^*)\half\EH(vv^*)\H(zz^*)\half\| \$\leq
  \|x^*\.y\|^2\ \|a^*\eH b\|^2\ \|\H(zz^*)\half\EH(ww^*)\H(zz^*)\half\| \$=
  \|x^*\.y\|^2\|a^*\eH b\|^2\|\H(zz^*)\half\EH(ww^*)\half\|^2 \$=
  \|x^*\.y\|^2\|a^*\eH b\|^2\|z^*\.w\|^2 =
  \|\xi\|^2 \|\phi\|^2 \|\eta\|^2,
  $$
  proving the statement when $N=1$.

  The method we will adopt to prove the general case will be to apply
the case we have already proved to matrices (see
\lcite{\MatrixInteraction}).  For arbitrary $N$ we will define
matrices $X, Y, Z, W, A, B \in M_n(A)$ to which we will apply the
single-summand matricial version of our result in order to obtain the
conclusion we seek.

The size of the matrices we will consider will be given by
$n=N^3$.  Also instead of using the index set $\{1,2,\ldots,N^3\}$
for columns and rows we will use the set
  $$
  \Lambda = \big\{ (i,j,k)\in\N^3 : 
    1\leq i\leq N,\ 
    1\leq j\leq N,\
    1\leq k\leq N\big\},
  $$
  which of course has $N^3$ elements.
  Therefore the matrix $X$, for example,  will be given in the form
  $$
  X = \Big(X_{(i,j,k)(l,m,n)}\Big)_{(i,j,k),(l,m,n)\in\hbox{$\Lambda$}}.
  $$
  For every $(i,j,k)$ and $(l,m,n)$ in  $\Lambda$ set
  $$
  X_{(i,j,k)(l,m,n)} \ =\ [1=l]\,[j=m]\,[k=n] \ x_i, \${}
  Y_{(i,j,k)(l,m,n)} \ =\ [1=l]\,[j=m]\,[k=n] \ y_i, \${}
  Z_{(i,j,k)(l,m,n)} \ =\ [i=l]\,[1=m]\,[k=n] \ z_j, \${}
  W_{(i,j,k)(l,m,n)} \ =\ [i=l]\,[1=m]\,[k=n] \ w_j, \${}
  A_{(i,j,k)(l,m,n)} \ =\ [i=l]\,[j=m]\,[1=n] \ a_k, \${}
  B_{(i,j,k)(l,m,n)} \ =\ [i=l]\,[j=m]\,[1=n] \ b_k,
  $$
  where the brackets refer to the obvious boolean valued function.

Observe that the coefficients above essentially mean that
$(i,j,k)=(l,m,n)$, except that in the first two $i$ is replaced by 1,
in the two middle ones $j$ is replaced by 1, and in the last two $k$
is replaced by 1.
  Therefore, for instance, $X_{(i,j,k)(l,m,n)}$ is almost always zero,
except when $l$ is equal to one, $j=m$, and $k=n$, in which case it is
equal to $x_i$.

We claim that the norm of $X^*\.Y$ (in the matricial version of $\X_0$) equals the norm of
  $\sum_{i} x_i^*\.y_i$ in $\X_0$.
  By definition
  $
  \|X^*\.Y\| = \|\H(XX^*)\half\EH(YY^*)\half\|,
  $ 
  and   we have
  $$
  \matrix{ (XX^*)_{(i,j,k)(l,m,n)} & = &
        \ds \sum_{(r,s,t)\in\Lambda} X_{(i,j,k)(r,s,t)}\big(X_{(l,m,n)(r,s,t)}\big)^* \hfill\cr\cr
  & = & \ds \sum_{(r,s,t)\in\Lambda} [1=r][j=s][k=t] [1=r][m=s][n=t] x_i  x_l^* \cr\cr
  & = & [j=m][k=n] x_i  x_l^*. \hfill}
  $$
  Likewise
  $
  (YY^*)_{(i,j,k)(l,m,n)} = [j=m][k=n] y_i  y_l^*.
  $
  Observe that the mapping
  $$\varphi: M_N(A) \to M_{N^3}(A)$$ which sends a matrix $u = (u_{ij})_{1\leq i,j\leq N}$ to
the matrix $\varphi(u)=\Big(U_{(i,j,k)(l,m,n)}\Big)_{(i,j,k),(l,m,n)\in\hbox{$\Lambda$}}$
given by
  $$
  U_{(i,j,k)(l,m,n)} = [j=m][k=n]u_{il},
  $$
  is an injective *-homomorphism and hence isometric. Moreover, by the
calculation just above, we have that $XX^*$ is the image of the
matrix $xx^* = (x_ix_j^*)_{1\leq i,j\leq N}$ 
under that mapping.  Clearly also $\varphi(\H(xx^*)) = \H(XX^*)$ and
$\varphi(\EH(yy^*)) = \EH(YY^*)$.
  So
  $$
  \|X^*\.Y\| =
  \|\H(XX^*)\half\EH(YY^*)\half\| =
  \|\varphi\big(\H(xx^*)\big)\half\varphi\big(\EH(yy^*)\big)\half\| \$=  
  \|\H(xx^*)\half\EH(yy^*)\half\| =  
  \Big\|\sum_{i} x_i^*\.y_i\Big\| =
  \|\xi\|,
  $$ 
  as desired.  In an entirely similar way one proves that
  $
  \|Z^*\.W\| =
  \|\eta\|,
  $
  and that 
  $
  \|A^*\eH B\| = \|\phi\|.
  $

  Applying the result proved in the first part to 
  $$
  \widetilde\xi  =  X^* \. Y, \quad
  \widetilde\eta =  Z^* \. W, \quad
  \widetilde\phi  =  A^* \eH B,
  $$
  we therefore conclude that
  $$
  \|\rp {\widetilde\xi}{R_{A^*,B}(\widetilde\eta)}\| \leq
  \|\widetilde\xi\|\, \|\widetilde\eta\|\, \|\widetilde\phi\| =
  \|\xi\|\, \|\eta\|\, \|\phi\|.
  \eqmark\AnInEquation
  $$  
  We now wish to compute $\rp
{\widetilde\xi}{R_{A^*,B}(\widetilde\eta)}$ explicitly.
  Letting $U$ and $V$ be defined by
  $U  =  \H(ZX^*)Y$, and $V = \EH(WA^*)B$, we have that
  $\rp{\widetilde\xi}{R_{A^*,B}(\widetilde\eta)} = U^*\eH V$, as in
\lcite{\Contas}.  We first compute
  $$
  U_{(l,m,n)(i,j,k)}  =
  \sum_{(p,q,r),(s,t,u)\in\Lambda}
    \H\big(Z_{(l,m,n)(p,q,r)}X_{(s,t,u)(p,q,r)}\null^*\big)Y_{(s,t,u)(i,j,k)} \$=
  \sum_{(p,q,r),(s,t,u)\in\Lambda}
    \tri(l,1,n)(p,q,r)
    \tri(1,t,u)(p,q,r)
    \tri(1,t,u)(i,j,k)
    \H(z_m x_s^*)y_s \$=
  [1=l][j=1][k=n][1=i]  \sum_{s=1}^N \H(z_m x_s^*)y_s.
  $$
  On the other hand
  $$
  V_{(l,m,n)(p,q,r)} =
  \sum_{(i,j,k),(s,t,u)\in\Lambda}
    \EH\big(W_{(l,m,n)(i,j,k)}A_{(s,t,u)(i,j,k)}\null^*\big)B_{(s,t,u)(p,q,r)} \$=
  \sum_{(i,j,k),(s,t,u)\in\Lambda}
    \tri(l,1,n)(i,j,k)
    \tri(s,t,1)(i,j,k)
    \tri(s,t,1)(p,q,r)
    \EH(w_ma_u^*)b_u \$=
  [p=l][q=1][1=n][1=r]  \sum_{u=1}^N
    \EH(w_ma_u^*)b_u.
  $$
  So  
  $$
  \Big(\rp{\widetilde\xi}{R_{A^*,B}(\widetilde\eta)}\Big)_{(i,j,k)(p,q,r)} = 
  \Big(U^*\eH V\Big)_{(i,j,k)(p,q,r)} \$=
  \sum_{(l,m,n)\in\Lambda} U_{(l,m,n)(i,j,k)}\null^* \eH V_{(l,m,n)(p,q,r)} \$=
  \sum_{(l,m,n)\in\Lambda}\ \sum_{s,u=1}^N\
      [1=l][j=1][k=n][1=i] 
      [p=l][q=1][1=n][1=r] $$ $$
      y_s^*\H(x_s z_m^*) \eH \EH(w_ma_u^*)b_u \$=
  \tri(i,j,k)(1,1,1)
  \tri(p,q,r)(1,1,1)
  \sum_{s,u,m=1}^N y_s^*\H(x_s z_m^*) \eH \EH(w_ma_u^*)b_u.
  $$
  This means that $\rp{\widetilde\xi}{R_{A^*,B}(\widetilde\eta)}$ is a matrix
with a single non-zero entry in the position $\big((1,1,1)(1,1,1)\big)$.
We claim that the value of this non-zero entry coincides with the term
in the left-hand-side of the inequation we are intent  on proving,
namely 
  $\RP \xi{\sum_{u=1}^N R_{a_u^*,b_u} (\eta)}.$
  In fact
  $$
  \RP \xi{\sum_{u=1}^N R_{a_u^*,b_u} (\eta)} =
  \RP {\sum_{s=1}^Nx_s^*\.y_s}{\sum_{u,m=1}^N R_{a_u^*,b_u} (z_m^*\.w_m)} \$=
  \sum_{s,u,m=1}^N \rp{x_s^*\.y_s}{z_m^*\V(w_ma_u^*)\.b_u} =
  \sum_{s,u,m=1}^N y_s^*\H\big(x_sz_m^*\V(w_ma_u^*)\big)\eH b_u,
  $$
  as desired.
  Therefore   
  $$
  \Big\|\RP \xi{\sum_{u=1}^N R_{a_u^*,b_u} (\eta)}\Big\| =
  \|\rp{\widetilde\xi}{R_{A^*,B}(\widetilde\eta)}\| 
  $$
  and the result follows from \lcite{\AnInEquation}.
  \proofend

  \state Corollary
  \label\TechCorollary
  Given $a_1,\ldots,a_N$ and $b_1,\ldots,b_N$ in $A$, let
$\phi=\sum_{k=1}^N a_k^* \eH b_k \in \KH$.  Then
  $$
  \Big\|\sum_{k=1}^N R_{a_k^*,b_k} (\eta)\Big\| \leq 
  \|\phi\|\,  \|\eta\|,
  $$
  for every $\eta\in\X_0$.

  \proof
  Letting  $\zeta=\sum_{k=1}^N R_{a_k^*,b_k} (\eta)$ we have that
  $$
  \|\zeta\|^2 =
  \|\rp \zeta\zeta\| =
  \Big\|\RP\zeta{\sum_{k=1}^N R_{a_k^*,b_k} (\eta)}\Big\| \leq
  \|\zeta\|\,  \|\eta\|\,  \|\phi\|,
  $$
  where the last inequality follows from \lcite{\TechLemma}.  So
  $
  \|\zeta\| \leq
  \|\eta\|\,  \|\phi\|,
  $
  as desired.
  \proofend

  Given $a,b\in A$ we therefore have that $R_{a,b}$
is a bounded operator on $\X_0$ and hence it extends to a bounded
operator (still denoted $R_{a,b}$) on $\X$, satisfying
  $$
  R_{a,b}(x\*y) = x\V(ya)\*b
  \for x,y\in A.
  $$

  Let $\D$ be the dense *-subalgebra of $\KH$ spanned by the set 
$\{a^*\eH b: a,b\in A\}$.  Each $\phi\in \D$ therefore  has the form
  $\phi  = \sum_{k=1}^N a_k^* \eH b_k$.  For each such $\phi$ define
an operator $R_\phi$ on $\X$ by
  $$
  R_\phi = \sum_{k=1}^N R_{a_k^*,b_k}.
  $$
  Clearly $R_\phi$ is a well-defined bounded operator on $\X$ but one
may suspect that the association $\phi \mapsto R_\phi$ is ill-defined.
Fortunately this is not so.    In fact observe that by 
\lcite{\TechCorollary} one has that 
  $
  \|R_\phi \| \leq 
  \|\phi\|,
  $
  and one may use this to show that $R_\phi$ is independent of the
given presentation of $\phi$.

  This shows not only that the association $\phi \mapsto R_\phi$ is
well-defined but also that it is continuous on the dense subalgebra
$\D$, therefore admitting an extension to the whole of $\KH$.

For $\phi\in \KH$
and $\eta\in \X$ we
will henceforth adopt the notation 
  $$
  \eta\phi := R_\phi(\eta).
  $$
  Since $\|R_\phi\|\leq\|\phi\|$ we have that
  $$
  \|\eta\phi\| \leq   \|\eta\|\,\|\phi\|
  \for \eta\in\X \for \phi\in\KH.
  $$
Observe that when $\eta=x\*y$, and $\phi=a\eH b$, we have that
$\eta\phi=R_{a,b}(\eta)$, so
  $$
  (x\*y)(a\eH b) =
  x\V(ya)\*b,
  $$
  which should be compared to \lcite{\SuggestRModStru}.

  \state Proposition
  With the above multiplication operation and the $\KH$--valued
inner-product $\rp{\,\cdot\,}{\,\cdot\,}$ described in
\lcite{\InnerProduct}, $\X$ is a \rhm{$\KH$}.
  
  \proof With respect to the module structure the only non-trivial
point remaining is the associativity axiom:
  $$
  (\eta\phi)\psi = \eta(\phi\psi)
  \for \eta\in\X\for \phi,\psi\in\KH.
  $$
  In order to prove this  we may invoke continuity to restrict ourselves to
the case in which  $\eta=x\*y$, $\phi=a\eH b$, and $\psi=c\eH d$, where
$x,y,a,b,c,d\in A$.
The right-hand side equals
  $$
  \eta(\phi\psi)=
  (x\*y)\big((a\eH b)(c\eH d)\big) =
  (x\*y)\big(a\EH(bc)\eH d\big) =
  x\V\big(ya\EH(bc)\big)\* d \$=
  x\V(ya)\V\big(\EH(bc)\big)\* d =
  x\V(ya)\V\H\V(bc)\* d =
  x\V(ya)\V(bc)\* d \$=
  (x\V(ya)\*b)(c\eH d) =
  \big((x\*y)(a\eH b)\big)(c\eH d) =
  (\eta\phi)\psi.
  $$
  We next prove that
  $$
  \rp \xi{\eta\phi} = \rp \xi\eta\phi
  \for\xi,\eta\in\X\for\phi\in\KH.
  $$
  Restricting, as we may, to the case in which $\xi=x\*y$,
$\eta=z\*w$, and $\phi=a\eH b$, where $x,y,z,w,a,b\in A$, we have
  $$
  \rp \xi{\eta\phi} =
  \rp {x\*y}{(z\*w)(a\eH b)} =
  \rp {x\*y}{z\V(wa)\* b} =
  y^*\H(x^*z\V(wa))\eH b \$=
  y^*\H(x^*z)\H\big(\V(wa)\big)\eH b =
  y^*\H(x^*z)\EH(wa)\eH b \$=
  \big(y^*\H(x^*z)\eH w\big)(a\eH b) =
  \rp{x\*y}{z\*w}(a\eH b) =
  \rp \xi\eta\phi.
  $$
  The positivity of the inner-product follows from \lcite{\Positivity}
and we leave it for the reader to verify the remaining details.
  \proofend

With respect to the
\lhm{$\KV$} structure of $\X$ we
restrict ourselves to the statement of the result, leaving it for the
reader to adapt the above reasoning to that situation.

  \state Proposition
  There is a unique
\lhm{$\KV$} structure on $\X$ with
respect to which
  $$
  (a\eV b)(x\*y) = a\*\H(bx)y
  \for a,b,x,y\in A.
  $$

  \state Proposition
  \label \ItIsABimodule
  With the above left and right module structures,  $\X$ is a
  \hbm{$\KV$}{$\KH$}.
  
  \proof All  we must prove  is that
  $$
  \lp\xi\eta\zeta =
  \xi\rp\eta\zeta
  \for \xi,\eta,\zeta\in\X,
  $$
  and we may clearly  restrict to the case in which $\xi=u\*v$,
$\eta=x\*y$, and $\zeta=z\*w$, where $u,v,x,y,z,w\in A$.  In this case
we have
  $$
  \xi\rp\eta\zeta =
  (u\*v)\rp{x\*y}{z\*w} =
  (u\*v)\big(y^*\H(x^*z)\eH w\big) =
  u\V\big(vy^*\H(x^*z)\big)\* w \$=
  u\V(vy^*)\V\big(\H(x^*z)\big)\* w \={(\Balanced.ii)}
  u\V(vy^*)\*\H(x^*z) w =
  \big(u\V(vy^*)\eV x^*\big)(z\*w) \$=
  \lp{u\*v}{x\*y}(z\*w) =
  \lp\xi\eta\zeta.
  \proofend
  $$

  We therefore reach the point where we may describe the ternary ring of
operators we are seeking:

  \state Proposition
  \label\TernaryOperation
  $\X$ is a ternary ring of operators with a unique ternary operation
  $$
  [\,\cdot\,,\,\cdot\,,\,\cdot\,]: \X\times\X\times\X \to \X
  $$
  such that 
  for every $u,v,x,y,z,w\in A$, one has  that
  $[u\*v, x\*y, z\*w] = u\V(vy^*)\*\H(x^*z) w.$

  \proof Any right-{\hm} is a ternary ring of operators with the operation 
  $$
  [x,y,z] := x\<y,z\>.
  \proofend
  $$

Whenever necessary  we will denote by $\X_\H$ the \rhm{$\KH$}
subjacent to $\X$ and similarly  $\X_\V$ will be the
corresponding \lhm{$\KV$}.

Recall that a {\hm} is said to be full when the ideal spanned by the
range of the inner-product is the whole coefficient algebra.  
  A quick glance at the definition of our inner-products (see
\lcite{\InnerProduct}) will reveal that both $\X_\V$ and $\X_\H$ are
full {\hm}s in view of the following:

  \state Lemma
  Let $E$ be a conditional expectation from a C*-algebra $A$ onto a
closed *-subalgebra $B$ and let $\K$ be the reduced C*-basic
construction with canonical idempotent $e$.  Then $\K$ is the closed
linear span of the set
  $$
  \{axeb: a,b\in A,\ x\in B\}.
  $$
    
  \proof
  Let $\{u_i\}$ be an approximate unit for $B$.  We claim that for
each $a,b\in A$ one has that
  $$
  aeb = \lim_{i\to\infty} au_ieb
  $$
  In fact
  $$
  \| aeb -au_ieb \|^2 =
  \| (a-au_i)eb \|^2 =
  \| b^*e(a^*-u_ia^*)(a-au_i)eb \|\$=
  \| b^*E(a^*a - a^*au_i - u_ia^*a + u_ia^*au_i)eb \| \$\leq
  \|E(a^*a) - E(a^*a)u_i - u_iE(a^*a) + u_iE(a^*a)u_i\|\, \|b\|^2,
  $$
  which converges to zero as $i\to\infty$.
  \proofend

En passant we observe that $\KV$
and $\KH$ are Morita-Rieffel equivalent \cite{\RieffelMorita}.

  This leads us to the following consequence of \scite{\BMS}{1.10}, in
which we denote by $\K(\M)$ the algebra of generalized compact
operators over a {\hm} $\M$:

  \state Proposition  
  \label\AntiAndIsomorphisms
  \izitem
  \zitem $\KV$ is isomorphic $\K(\X_H)$.
  \zitem $\KH$ is anti-isomorphic to  $\K(\X_V)$.

In  \cite{\Kasparov} it is proved that the algebra of all adjointable
operators on a {\hm} is isomorphic the the multiplier algebra
of the algebra of generalized compact operators.  Therefore, denoting
algebras of adjointable operators by $\L(\,\cdot\,)$ and 
multiplier algebras by $\M(\,\cdot\,)$, we have by
\lcite{\AntiAndIsomorphisms} that
  $$
  \M(\KV) \sim \M(\K(\X_H)) \sim \L(\X_H),
  $$
  while
  $$
  \M(\KH) \antiiso \M(\K(\X_V)) \sim \L(\X_V),
  $$
  where $\antiiso$ means anti-isomorphic.  One may then consider $\X$
as a \hbm{$\M(\KV)$}{$\M(\KH)$}.  

Using  the remark after Lemma (2.1.3) in \cite{\Watatani} it is easy
to see that there are
*-homomorphisms
  $$
  \lambda_\V : A \to \M(\KV)
  \and
  \lambda_\H : A \to \M(\KH),
  $$
  such that
  $$
  \lambda_\V(a)x\eV y = ax\eV y
  \and
  \lambda_\H(a)x\eH y = ax\eH y,
  $$
  for all $a,x,y\in A$.  
  We may therefore view $\X$ as an $A$--$A$--bimodule, although not
necessarily as a {\hb} since there is no guarantee that the
inner-products lie within $\lambda_\V(A)$ of $\lambda_\H(A)$.

  In order to obtain concrete formulas for the left and right action
of $A$ on $\X$ observe that for all  $a,b,c,x,y\in A$ one has 
  $$
  \big((b\eV c)\lambda_\V(a)\big)(x\*y) =  
  (b\eV ca)(x\*y) =  
  b\*\H(cax)y =  
  (b\eV c)(ax\*y).
  $$
  So, for all $k\in\KV$, we have $k\lambda_\V(a)(x\*y) = k(ax\*y)$ and
hence $\lambda_\V(a)(x\*y)=ax\*y$.  In other words, the left action of
$A$ on $\X$ is given simply by
  $$
  a\cdot (x\*y) = (ax)\*y
  \for a,x,y\in A.
  $$
  Likewise the right action of $A$ on $\X$ is given by
  $$
  (x\*y)\cdot a  = x\*(ya)
  \for a,x,y\in A.
  $$

  The relationship between the TRO structure of $\X$ and its
\bm{$A$}{$A$} is given by our next:

  \state Proposition
  \label\ExampleCorrespondence
  For every $a\in A$ and every $\xi,\eta,\zeta\in\X$ one has that
  \izitem
  \zitem $[\xi,a\eta,\zeta] =   [\xi,\eta,a^*\zeta]$, and
  \zitem $[\xi,\eta a,\zeta] =   [\xi a^*,\eta,\zeta]$.

  \proof Immediate from the description of the ternary product in
\lcite{\TernaryOperation}.
  \proofend

  \section{Constructing a covariant representation}
  As in the previous section we fix an interaction $(\V,\H)$ over a
C*-algebra $A$.   However, for simplicity, we now assume that $A$ is unital.

 We are now ready to achieve the construction of a nondegenerated
covariant representation for $(\V,\H)$

  Considering the \hbm{$\M(\KV)$}{$\M(\KH)$} structure of $\X$
described in section \lcite{\ReconstructionSection}, we will let
$\Link$ be its linking algebra \scite{\BMS}{2.2}.  Thus $\Link$ may be
represented as
  $$
  \Link = 
  \left[\matrix{\M(\KV) & \X \cr \X^* & \M(\KH) \cr}\right].
  $$

Since $A$ is assumed to be unital one may speak of the element
$1\*1\in\X$.
  Speaking of 1 we will shortly need the following simple technical
fact:

  \state Lemma
  \label\LemaAboutOne
  If $A$ is unital then $\V(1)\eV=\eV$ in $\KV$, and $\H(1)\eH=\eH$ in
$\KH$.

  \proof Recalling the construction of $\KV$ as operators on the \rhm{$\V(A)$}
obtained by equipping $A$  with the inner-product 
  $$
  \<a,b\> = \EV(a^*b)
  \for a,b\in A,
  $$
  we have for all $a\in A$ that
  $$
  \Big(\V(1)\eV\Big)\calcat a =
  \V(1)\EV(a) =
  \V(1)\V\big(\H(a)\big) \={(\DefInteraction.iv)}
  \V\big(1\H(a)\big) =
  \EV(a) =
  \eV\calcat a,
  $$
  proving that 
  $\V(1)\eV=\eV$.  That  $\H(1)\eH=\eH$ follows by a similar argument.
  \proofend

  \state Proposition
  \label \MainResult
  Assume that $A$ is unital and 
  let
  $
  \pi : A \to \Link
  $
  be the  *-homomorphism 
  given by
  $$
  \pi(a) =
  \left[\matrix{\lambda_\V(a) & 0 \cr 0 & \lambda_\H(a) \cr}\right]
  \for a\in A,
  $$
  and 
  let $\S$ be the element of $\Link$ given by
  $$
  \S=  
  \left[\matrix{0 & 1\*1 \cr 0 & 0 \cr}\right].
  $$
  Then $(\pi,\S,\Link)$ is a nondegenerated covariant representation of $(\V,\H)$.

  \proof By direct computation we conclude that $\S\pi(a)\S^*$ is a matrix
having a single nonzero entry in the top left-hand corner and the value of
that entry equals
  $$
  (1\*1)\lambda_\H(a)(1\*1)^* = 
  (1\*a)(1\*1)^* = 
  \lp{1\*a}{1\*1} =
  \V(a)\eV.
  $$
  The same holds for $\pi(\V(A))\S\S^*$, except that the
single nonzero entry in the top left-hand corner is
  $$
  \lambda_\V\big(\V(a)\big)(1\*1)(1\*1)^* =
  \lambda_\V\big(\V(a)\big)\V(1)\eV =
  \V(a)\V(1)\eV \={(\LemaAboutOne)}
  \V(a)\eV.
  $$
  This shows \lcite{\DefCovarRep.i} and the proof of
\lcite{\DefCovarRep.ii} may be similarly given.  

Let us now show nondegeneracy.  Focusing on \lcite{\DefineNonDeg.i}
let $a\in\V(A)$ be such that $\pi(a)\S=0$.  Then $\pi(a)\S\S^*=0$, as
well, but direct computation shows that this is a matrix whose top
left-hand corner is given by $a\eV$.  But then
\scite{\Watatani}{2.1.4(2)} implies that $a=0$.
  \proofend

Supposing that $\EV$ or $\EH$ are faithful one then has that $\pi$ is
faithful by the argument presented shortly after \lcite{\DefineNonDeg}.
Therefore $A$ may be identified with its image in $\Link$ which, together
with $\S$, satisfies \lcite{\CovarianceAxiom}.

However we do not needed faithfulness in order to prove:

  \state Theorem
  Let $A$ be a unital C*-algebra and let $(\V,\H)$ be an interaction
over $A$.  Then there exists a C*-algebra $B$ containing $A$ and a
partial isometry $\S\in B$ such that axioms \lcite{\CovarianceAxiom} and
\lcite{\NonDegeneracyAxiom} hold.

  \proof
  Let $(\pi,\S,B)$ be a nondegenerated covariant representation of
$(\V,\H)$, such as the one constructed in \lcite{\MainResult}
Now consider the covariant representation
  $$
  (\pi\oplus id,\S\oplus0,B\oplus A)
  $$
  where $id$ is the identity map on $A$.  It is elementary to verify
that this is indeed a covariant representation.  Evidently $\pi\oplus
id$ is injective and hence we may identify $A$ with its image in
$B\oplus A$.  With the identification 
\lcite{\CovarianceAxiom} become a consequence of the properties of a
covariant representation.  With respect to \lcite{\NonDegeneracyAxiom.i}
observe that, if $a\in\V(a)$ and 
$(\pi(a)\oplus a)(\S\oplus 0)=0$, then
$\pi(a)S = 0$ and hence $a=0$ by nondegeneracy of $\pi$.  Likewise one
shows \lcite{\NonDegeneracyAxiom.ii}.
   \proofend

The reader may be struck with the impression that the wild juggling of
covariant representations (consider $\S=0$ for instance!) is a bit
exaggerated and that something must be done to stop it.  I agree.  But
so far I don't know exactly how.  However I have a hunch as to where
to look for a solution.

  \section{{\Mycorr}s and the covariance algebra}
  In this section we wish to discuss a construction which should be
thought of as generalizing the crossed-product in case the given
interaction is a pair $(\a,\Tr)$ consisting of an endomorphism $\a$
and a transfer operator $\Tr$ as in \cite{\endo}.  Since this can be
done in the much larger generality of {\mycorr}s we will broaden our
context to include the latter situation.

Although we hope this construction will catch the attention of a
number of readers, we still consider it as an experimental attempt,
especially given that results presented here will mostly be of a formal
nature.  

  \definition
  \label \DefSymCorr
  Let $A$ be a C*-algebra.  A \stress{\mycorr} over $A$ is a ternary
ring of operators $\X$, with ternary operation
$[\,\cdot\,,\,\cdot\,,\,\cdot\,]$, which is at the same time an
$A$--$A$--bimodule such that for every $a\in A$, and
$\xi,\eta,\zeta\in\X$ one has that
  \izitem
  \zitem $[\xi,a\eta,\zeta] =   [\xi,\eta,a^*\zeta]$, and
  \zitem $[\xi,\eta a,\zeta] =   [\xi a^*,\eta,\zeta]$.
  
Observe that the TRO constructed in section
\lcite{\ReconstructionSection} out of an interaction is an example of
this situation.

\medskip Given a {\mycorr} $\X$ over $A$ observe that, using the
terminology of \scite{\TPA}{4.3}, for each $a\in A$  one has that the operator of left
multiplication by  $a$, namely the operator
  $$
  \lambda(a) : \xi \in \X  \mapsto a\xi \in\X
  $$
  is a \stress{left operator}, while
  $$
  \rho(a) : \xi \in \X  \mapsto \xi a \in\X
  $$
  is a \stress{right operator}.

If $X$ is any TRO whatsoever we shall denote by $\L(\X)$ the C*-algebra \scite{\TPA}{4.5}
of all left operators on $\X$ and by $\R(\X)$ the C*-algebra of all
right operators.

  \state Lemma If $\X$ is a TRO then every right operator  commutes with
every left operator.

  \proof  Let $R$ be a right operator  and $L$ be a left operator.
Then, for every $\xi,\eta,\zeta\in\X$ one has that
  $$
  [\xi,R(L(\eta)),\zeta] =
  [R^*(\xi),L(\eta),\zeta] =
  [R^*(\xi),\eta,L^*(\zeta)] =
  [\xi,R(\eta),L^*(\zeta)] =
  [\xi,L(R(\eta)),\zeta].
  $$
  Thus, if $\delta=R(L(\eta))-L(R(\eta))$, one has that  $[\xi,\delta,\zeta]=0$, for
every $\xi,\zeta\in\X$.  In particular
  $$
  \|\delta\|^3 = \|[\delta,\delta,\delta]\| = 0.
  \proofend
  $$
  
It follows that a TRO $\X$ is always an \bm{$\L(\X)$}{$\R(\X)^{op}$} (the
superscript ``$^{op}$'' standing for the opposite C*-algebra).

  \definition Let $\X$ be a TRO.  Given $\xi,\eta\in \X$, we will
denote by $\theta^r_{\xi,\eta}$ and $\theta^\ell_{\xi,\eta}$ the operators on $\X$ defined by
  $$
  \theta^r_{\xi,\eta} : x \in \X \mapsto [x,\xi,\eta]\in\X,
  $$ and $$
  \theta^\ell_{\xi,\eta} : x \in \X \mapsto [\xi,\eta,x]\in\X.
  $$

By the associativity property of TRO's we have that
  $$
  [x,\theta^\ell_{\xi,\eta}(y),z] = 
  [x,[\xi,\eta,y],z] = 
  [x,y,[\eta,\xi,z]] = 
  [x,y,\theta^\ell_{\eta,\xi}(z)],
  $$
  and hence $\theta^\ell_{\xi,\eta}$ is a left operator with
$(\theta^\ell_{\xi,\eta})^*=\theta^\ell_{\eta,\xi}$.
  Likewise $\theta^r_{\eta,\xi}$ is a right operator and
$(\theta^r_{\xi,\eta})^*=\theta^r_{\eta,\xi}$.

\definition  We will denote by $\KL(\X)$ the closed linear span, within $\L(\X)$,
of the set $\{\theta^\ell_{\xi,\eta}: \xi,\eta \in\X\}$, and by $\KR(\X)$
the closed linear span, within $\R(\X)$, of $\{\theta^r_{\xi,\eta}:
\xi,\eta \in\X\}$.  The operators in $\KL(\X)$ will be called
\stress{generalized compact left operators} and those in $\KR(\X)$ will be
called \stress{generalized compact right operators}.

  Defining an $\L(\X)$--valued inner-product
$\lp{\,\cdot\,}{\,\cdot\,}$ on $\X$ by 
  $$
  \lp\xi\eta=
  \theta^\ell_{\xi,\eta}
  \for \xi,\eta\in\X,
  $$
  and an  $\R(\X)$--valued inner-product
$\rp{\,\cdot\,}{\,\cdot\,}$ on $\X$ by 
  $$
  \lp\xi\eta=
  \theta^r_{\xi,\eta}
  \for \xi,\eta\in\X,
  $$
  we have that $\X$ is an \hbm{$\L(\X)$}{$\R(\X)^{op}$}.  The crucial
positivity axiom holds here precisely because $\X$ is assumed to be a
TRO, as opposed to a general ternary C*-ring.  See \cite{\Zettl} or
\cite{\TPA} for more details.

  Now let $\X$ be a {\mycorr} over $A$.  The \bm{$A$}{$A$} structure
of $\X$ therefore gives *-homomorphisms
  $$ 
  \lambda : A \to \L(\X)
  \and
  \rho: A \to \R(\X)^{op}.
  $$

Conversely, suppose that $L$ and $R$ are C*-algebras and $\X$ is an
\hbm{$L$}{$R$} with $L$--valued inner-product
$\lp{\,\cdot\,}{\,\cdot\,}$,
and $R$--valued inner-product $\rp{\,\cdot\,}{\,\cdot\,}$.
Suppose moreover that we are given *-homomor\-phisms
  $$
  \lambda: A \to L
  \and
  \rho: A \to R.
  $$
Consider $\X$ as an $A$--$A$--bimodule via $\lambda$ and $\rho$, and
as a TRO under the ternary operation
  $$
  [\xi,\eta,\zeta] = \xi\rp\eta\zeta
  \for \xi,\eta,\zeta\in\X.
  \eqno{(\seqnumbering)}
  \label \TROHBOp
  $$

We claim that this makes $\X$ into a {\mycorr}.  In fact, in order to
check that \lcite{\DefSymCorr.i} holds, let $a\in A$ and
$\xi,\eta,\zeta\in\X$.  Then, recalling that $\xi\rp\eta\zeta =
\lp\xi\eta\zeta$, we have
  $$
  [\xi,a\eta,\zeta] = 
  \xi\rp{a\eta}\zeta =
  \lp\xi{a\eta}\zeta =
  \lp\xi\eta a^*\zeta =
  \xi\rp\eta{a^*\zeta} =
  [\xi,\eta,a^*\zeta].
  $$
  As for \lcite{\DefSymCorr.ii},
  $$
  [\xi,\eta a,\zeta] = 
  \xi\rp{\eta a}\zeta =
  \xi a^*\rp\eta\zeta =
  [\xi a^*,\eta,\zeta].
  $$

We therefore obtain the following alternative characterization of {\mycorr}s:

  \state Proposition
  $\X$ is a {\mycorr} over $A$ if and only if $\X$ is
an \hbm{$L$}{$R$} for C*-algebras $L$ and $R$, and there are
*-homomorphisms $\lambda: A \to L$ and $\rho: A \to R$ such that
the \bm{$A$}{$A$} structure of $\X$ is given by $\lambda$ and $\rho$
and the ternary operation on $\X$ is given as in \lcite{\TROHBOp}.

  \definition Let $\X$ be a {\mycorr} over the C*-algebra $A$.  A
\stress{covariant representation of the pair $(A,\X)$ in a C*-algebra $B$} is a
pair $(\pi,\psi)$, where $\pi:A\to B$ is a *-homomorphism and $\psi:\X\to
B$ is a contractive linear map, such that for all $a\in A$, and
$\xi,\eta,\zeta\in\X$ one has that
  \izitem
  \zitem $\pi(a)\psi(\xi) = \psi(a \xi)$,
  \zitem $\psi(\xi)\pi(a) = \psi(\xi a)$,
  \zitem $\psi([\xi,\eta,\zeta]) = \psi(\xi)\psi(\eta)^*\psi(\zeta)$.

  \medskip Fix, for the time being, a {\mycorr} $\X$ over the
C*-algebra $A$.
  By \cite{\Blackadar} we therefore have that there exists a
universal C*-algebra for covariant representations.  This algebra will
be denoted by $\TAX$ and it will be referred to as the Toeplitz algebra
for the pair $(A,\X)$.  It is the codomain of a canonical covariant
representation $(\pi_u,\psi_u)$ such that for every C*-algebra $B$ and
any covariant representation $(\pi,\psi)$ in $B$ there exists a unique
*-homomorphism $\phi:\TAX \to B$ such that for every $a\in A$ and
$\xi\in\X$ one has that
  \izitem   
  \zitem $\phi(\pi_u(a)) = \pi(a)$, and
  \zitem $\phi(\psi_u(\xi)) = \psi(\xi)$.

  \medskip Here is a somewhat silly example of covariant
representation: take $\pi$ to be any *-homomorphism from $A$ to a
C*-algebra $B$ and set $\psi=0$.  This at least has the merit of
showing that $\pi_u$ is injective, since $\pi_u(a)=0$ would imply that
$\pi(a)=0$ and, supposing $\pi$  injective, we have that $a=0$.

  A more interesting example of covariant representation is suggested
by \lcite{\MainResult}: view $\X$ as an \hbm{$\L(\X)$}{$\R(\X)^{op}$}
as above and let $\Link$ be the linking algebra of $\X$.  Define
  $$
  \psi:\xi\in\X\mapsto 
  \left[\matrix{0 & \xi \cr 0 & 0}\right] \in \Link.
  $$
  Also let $\lambda:A\to\L(\X)$ and $\rho:A\to\R(\X)$ be the maps corresponding
to the \bm{$A$}{$A$} structure of $\X$ and set 
  $$
  \pi: a\in A \mapsto
  \left[\matrix{\lambda(a) &0 \cr 0 & \rho(a)}\right] \in \Link.
  $$
  It is then elementary to check that $(\pi,\psi)$ is a covariant
representation. 
Since $\psi$ is isometric it follows that $\psi_u$ is
isometric as well.

Since $\pi_u$ and $\psi_u$ are isometric we may identify  both $A$ and
$\X$ with their images within $\TAX$ via $\pi_u$ and $\psi_u$, respectively.
The ternary operation of $\X$, seen within $\TAX$, is therefore given by
  $$
  [\xi,\eta,\zeta] =   \xi\eta^*\zeta
  \for \xi,\eta,\zeta\in\X.
  $$
  Given $\xi_1,\xi_2,\ldots,\xi_n$ and $\eta_1,\eta_2,\ldots,\eta_n$
in $\X$, let $k$ be the generalized compact left operator given by 
$k = \sum_{i=1}^n\theta^\ell_{\xi_i,\eta_i}$.  For every $x\in\X$ we then
have
  $$
  k(x) = 
  \sum_{i=1}^n\theta^\ell_{\xi_i,\eta_i}(x) =
  \sum_{i=1}^n[\xi_i,\eta_i,x] =
  \Big(\sum_{i=1}^n\xi_i\eta_i^*\Big)x,
  $$
  and hence the action of $k$ on $\X$ is identical to the operator of
left multiplication by   $\sum_{i=1}^n\xi_i\eta_i^*$.  By
\scite{\TPA}{4.7} we have that 
  $$
  \|k\| = \Big\|\sum_{i=1}^n\xi_i\eta_i^*\Big\|,
  $$
  which implies that there exists an isometric linear map
$\Lambda:\KL(\X)\to \TAX$ such that
$\Lambda(\theta^\ell_{\xi,\eta})=\xi\eta^*$.  It is elementary to show
that this is a *-homomorphism as well and hence we may identify
$\KL(\X)$ with the subset of $\TAX$ given by $\X\X^*$ (closed linear
span).  Similarly we will identify $\KR(\X)$ with $\X^*\X$.

One of the reasons for introducing the notion of {\mycorr} is that it
generalizes the well known notion of correspondences \cite{\Muhly},
which in turn was motivated by the so called Pimsner bimodules
\cite{\Pimsner}.  In fact, let $\X$ be a correspondence over $A$, that
is, $\X$ is a \rhm{$A$} equipped with a *-homomorphism $\lambda$ from
$A$ into the C*-algebra $\L_A(\X)$ of all adjointable operators on
$X$.  Then $\X$ is clearly an \hbm{$A$}{$A$} and a TRO with the
ternary operation defined in \lcite{\TROHBOp}, and it is elementary to
check that \lcite{\DefSymCorr.i-ii} hold.  Summarizing $\X$ is a
{\mycorr}.

  In this case notice that for all $\xi,\eta\in \X$, the
generalized right compact operator $\theta^r_{\xi,\eta}$ lies in the
range of $\rho$.  In fact, for all $x\in\X$ we have
  $$
  \theta^r_{\xi,\eta}(x) = 
  [x,\xi,\eta] =
  x\rp\xi\eta =
  \rho\big(\rp\xi\eta\big)\calcat x,
  $$
  so that 
  $$
  \theta^r_{\xi,\eta} = \rho\big(\rp\xi\eta\big).
  \eqno{(\seqnumbering)} 
  \label \KRinA
  $$
  Since the
range of a *-homomorphism is always closed, we conclude that
  $
  \KR(\X)\subseteq\rho(A).
  $

  We feel that the most significant difference between the usual
notion of correspondence and its generalized version is that the
above inclusion is not required to hold in general.  But even in the
general case one should pay attention to the relationship between
generalized right compact operators and $\rho(A)$:

  \definition Let $\X$ be a {\mycorr} over the C*-algebra $A$.  By a
\stress{right redundancy} we shall mean a pair
$(a,k)\in\TAX\times\TAX$, such that $a\in A$, $k\in\X^*\X\
(=\KR(\X))$, and
  $$
  xa = xk \ (= k(x))
  \for x\in\X.
  $$
  Likewise, by a \stress{left redundancy} we shall mean a pair
$(a,k)\in\TAX\times\TAX$, such that $a\in A$, $k\in\X\X^*\ (=\KL(\X))$, and
  $$
   ax  =  kx\ (= k(x))
  \for x\in \X.
  $$

  Observe that the right redundancies are precisely the pairs of the
form $(a,\rho(a))$, where $a\in\rho\inv(\KR(\X))$.  Symmetrically the
left redundancies are the pairs $(a,\lambda(a))$, where
$a\in\lambda\inv(\KL(\X))$.
  For the case of left redundancies this is related to condition (4)
of \scite{\Pimsner}{Definition 3.8}.

  \definition The \stress{left} (resp.~\stress{right})
\stress{redundancy ideal}, denoted $\J_\ell$ (resp.~$\J_r$), is the
closed two-sided ideal of $\TAX$ generated by the elements $a-k$, for
every left (resp.~right) redundancy $(a,k)$ such that $a$ lies in
$\Ker(\lambda)^\perp$ (resp.~$\Ker(\rho)^\perp$), the symbol $^\perp$
standing for the annihilator of the given ideal.  The \stress{full
redundancy ideal} $\J$ is, by definition, the closed two-sided ideal
generated by $\J_\ell\cup\J_r$.

Restricting the element $a$ above to the annihilator of, say
$\Ker(\lambda)$, is an idea of Katsura \scite{\Katsura}{2.3, 2.5}
and it has the purpose of avoiding an otherwise nonzero intersection
between $A$ and $\J^\ell$ (which we will later see is undesirable).
In fact, let $(a,k)$ be a left redundancy such that
$a\notin\Ker(\lambda)^\perp$ and suppose that $a-k$ belong to a given
ideal $\I$ of $\TAX$.  Then for some $b\in\Ker(\lambda)$ we have that
$ba\ne0$ and for every $\xi\in\X$ we have
  $$
  0 = \lambda(ba)\xi = ba\xi = bk(\xi).
  $$
  So $bk=0$ and thus $ba = b(a-k)\in\I.$ Therefore $0\neq ba\in\I\cap
A$ and hence $\I\cap A\ne\{0\}$.
  Fortunately, in case both $\lambda$ and $\rho$ are injective, we do
not have to worry about this point.

Let us now discuss a relationship between our construction and the
well known construction of Cuntz-Pimsner algebras.  For this let $A$
be a C*-algebra and let $\X$ be a correspondence in the usual sense.
For simplicity we will assume that the left action $\lambda$ is
injective and the {right-\hm} is full in the sense that  $A$ coincides
with the closed linear span of $\rp\X\X$ (in which case $\rho$ must be
injective%
  \fn{Suppose $a$ is such that $\rho(a)=0$. Then for every
  $\xi,\eta\in\Xs$ we have that
  $
  0 =
  \rp\xi{\rho(a)\eta} =
  \rp\xi{\eta a} =
  \rp\xi\eta a,
  $
  and hence $0=\rp{\Xs}{\Xs\,} a=Aa$,
  which implies that $a=0$.  See also \scite{\BMS}{1.10}.
  }
  as well).  These are the hypotheses under which Pimsner introduced
the construction in \cite{\Pimsner} of what is now known as the
Cuntz-Pimsner algebras (see the first paragraph in
\scite{\Pimsner}{Chapter 1} as well as \scite{\Pimsner}{Remark 1.2.3}).

Let us denote by $\TX$ and by $\OX$ the Toeplitz-Cuntz-Pimsner and the
Cuntz-Pimsner algebras associated to $\X$, respectively.

Seeing, as before, $\X$ as a {\mycorr},  our aim is to show:

  \state  Proposition
  \label \CuntzPimsnerAndMe
  Under the above assumptions $\TX$ is isomorphic to the quotient of\/
$\TAX$ by the right redundancy ideal and $\OX$ is the quotient of\/
$\TAX$ by the full redundancy ideal.

  \proof By the universal property of $\TAX$ there exists a
*-homomorphism 
  $$
  f:\TAX\to\TX
  $$
  such that $f(z)=z$, for all $z\in A\cup\X$.

Given an element $a\in A$ of the form $a=\rp\xi\eta$ we have seen in
\lcite{\KRinA} that $\theta^r_{\xi,\eta} = \rho\big(\rp\xi\eta\big)$.
The assumption that $A$ is spanned by $\rp\X\X$ therefore implies that
$\KR(\X)=\rho(A)$ and hence the right redundancies are precisely the
pairs $(a,\rho(a))$ for $a\in A$.  The right redundancy ideal is
therefore spanned by the differences $\rp\xi\eta-\xi^*\eta$, for
$\xi,\eta\in\X$.  We then have that $f$ vanishes on $\J_r$ by the very
definition of $\TX$ (see \scite{\Pimsner}{3.4.(3)}).  By passage to the
quotient this gives a *-homomorphism
  $\widetilde f:\TAX/\J_r\to\TX$ which we will prove to be a
bijection.

Composing the canonical maps $\pi_u$ and $\psi_u$ with the quotient
map we get
  $$
  \widetilde\pi_u: A \to \TAX/\J_r
  \and 
  \widetilde\psi_u: \X \to \TAX/\J.
  $$
  Again because $\rp\xi\eta-\xi^*\eta$ is in $\J_r$ we conclude that
$\widetilde\pi_u(\rp\xi\eta) =
\widetilde\psi_u(\xi)^*\widetilde\psi_u(\eta)$ and hence the pair
$(\widetilde\pi_u, \widetilde\psi_u)$ satisfies the hypotheses in
\scite{\Pimsner}{3.4} thus giving a *-homomorphism $g:\TX\to\TAX$ which
is easily seen to be the inverse of $\widetilde f$.

The quotient $\TAX$ by the full redundancy ideal is clearly isomorphic
to the quotient of $\TAX/\J_r$ by the quotient image of $\J_\ell$.
Equivalently it is isomorphic  to the quotient of $\TX$ by the corresponding ideal.
This is well known to be isomorphic to $\OX$.
  \proofend

This maybe justifies giving the following:

  \definition Let $\X$ be a {\mycorr} over a C*-algebra $A$.  Then the
covariance algebra for the pair $(A,\X)$, denoted $C^*(A,\X)$, is the
quotient of $\TAX$ by the full redundancy ideal.

As already mentioned this is a very tentative definition which we
nevertheless believe to be of interest given its similarities with the
enormously popular Cuntz-Pimsner construction.  The main questions
risen by the above definition are:

\medskip (1) Are the canonical embeddings of $A$ and $\X$ into
$C^*(A,\X)$ injective?  Equivalently, is the intersection between the
redundancy ideal $\J$ and $A$, or $\X$, the zero ideal?

\medskip  (2) Can one say anything useful about the structure of\/ $C^*(A,\X)$,
compute its $K$-theory, or find a concrete faithful representation of it?

\medskip (3) Is there a \stress{Fock space representation} of
$\T(A,\X)$ similar to the one given in \cite{\Pimsner}?

\medskip As far as examples are concerned let $(\V,\H)$ be an
interaction over a C*-algebra $A$.  Then by \lcite{\TernaryOperation}
and \lcite{\ExampleCorrespondence} we have that the ternary ring of
operators $\X = \XVH$ constructed in section
\lcite{\ReconstructionSection} is a {\mycorr}.  One could then attempt
to study the \stress{crossed-product algebra} $C^*(A,\V,\H)$, namely
the covariance algebra $C^*(A,\XVH)$.  This at least includes
our previous crossed-product construction:

  \state Proposition
  Let $A$ be a unital C*-algebra, $\alpha$ be an injective
endomorphism of $A$, and $\Tr$ be a 
  nondegenerated\fn{Meaning that $\Tr(a^*a)=0$ implies that $a=0$.}
  transfer operator with $\Tr(1)=1$.  Then, viewing $(\a,\Tr)$ as an
interaction over $A$ by \lcite{\TransferInteraction} one has that the
covariance algebra $C^*(A,\XAT)$ coincides with the crossed-product
$A\crossproduct_{\a,\Tr}\N$.

  \proof
  For all $a,b\in A$ notice that
  $$
  a\*b =
  a\*\Tr(\a(b))1 \={(\Balanced)}
  a\a(b)\*1.
  $$
  So $A\*1$ is dense in $\XAT$.  Moreover
  by \lcite{\NormComputationAbstract}
  one has for all $a\in A$ that
  $$
  \|a\*1\| =
  \|\Tr(a^*a)\|\half.
  $$
  So it follows that the map
  $$
  \phi: m\in A \mapsto m\*1\in A\.A
  $$
  extends to an isometric linear map from the \rhm{$A$} $\M_\Tr$,
defined in \scite{\endo}{Section 3}, onto $\XAT$.  We claim that this
is an isomorphism of {\mycorr}s, meaning that it is $A$--$A$-linear
and preserves the ternary operation.
  In fact, for $m\in A\subseteq\M_\Tr$ and $a\in A$ we have
  $$
  \phi(m\cdot a) = 
  \phi(m\a(a)) =
  m\a(a)\*1 =
  m\*a =
  (m\*1)a =
  \phi(m)a,
  $$
  so $\phi$ is right--$A$--linear.  Also
  $$
  \phi(a\cdot m) = 
  \phi(a m) = 
  a m\*1 =
  a(m\*1) =
  a\phi(m),
  $$  
  establishing the left--$A$--linearity.  As for the ternary operation
we have for $x,y,z\in A\subseteq\M_\Tr$ that
  $$
  \phi\big([x,y,z]\big) =
  \phi\big(x\cdot\rp yz\big) =
  \phi\big(x\cdot\Tr(y^*z)\big) =
  \phi\big(x\a(\Tr(y^*z))\big) \$=
  x\a(\Tr(y^*z))\*1 \={(\Balanced.ii)}
  x\*\Tr(y^*z),
  $$
  while
  $$
  [\phi(x),\phi(y),\phi(z)]=
  [x\*1,y\*1,z\*1] =
  x\a(11^*)\*\Tr(y^*z) 1 =
  x\*\Tr(y^*z),
  $$
  proving the claim.  Therefore
  \def\iso{\ \simeq\ }
  $$
  A\crossproduct_{\a,\Tr}\N \iso
  \O_{\M_{\Tr}} \iso
  \O_{\Xs_{\a,\Tr}} {\buildrel {(\CuntzPimsnerAndMe)} \over {\iso}}
  C^*(A,\XAT).
  $$
  For the first isomorphism above see also \cite{\Raeburn}.
  \proofend

\endgroup

  \bigskip \bigskip
  \font \sc = cmcsc8 \sc

  \noindent Departamento de Matem\'atica,
  Universidade Federal de Santa Catarina

  \noindent  88040-900 --  Florian\'opolis --
  Brasil (exel@mtm.ufsc.br).

\end
\references

\bibitem{\Blackadar}
  {B. Blackadar}
  {Shape theory for $C^*$-algebras}
  {{\it Math. Scand.}, {\bf 56} (1985), 249--275}

\bibitem{\BMS}
  {L. G. Brown, J. A. Mingo, and N. T. Shen}
  {Quasi-multipliers and embeddings of Hilbert $C^*$-bimodules}
  {{\it Canad. J. Math.}, {\bf 46} (1994), 1150--1174}

\bibitem{\Raeburn}
  {N. Brownlowe and I. Raeburn}
  {Exel's Crossed Product and Relative Cuntz-Pimsner Algebras}
  {preprint, 2004, [arXiv:math.OA/0408324]}

\bibitem{\Deaconu}
  {V. Deaconu}
  {C*-algebras of commuting endomorphisms}
  {preprint, 2004, [arXiv:math.OA/0406624]}

\bibitem{\TPA}
  {R. Exel}
  {Twisted partial actions, a classification of regular C*-algebraic bundles}
  {{\it Proc. London Math. Soc.}, {\bf 74} (1997), 417--443,
[arXiv:funct-an/9405001], MR 98d:46075}

\bibitem{\endo}
  {R. Exel}
  {A new look at the crossed-product of a C*-algebra by an endomorphism}
  {{\it Ergodic Theory Dynam. Systems}, to appear, [arXiv:math.OA/0012084]}

\bibitem{\tower}
  {R. Exel}
  {Crossed-products by finite index endomorphisms and KMS states}
  {{\it J. Funct. Analysis}, {\bf 199} (2003), 153--188, [arXiv:math.OA/0105195]}

\bibitem{\vershik}
  {R. Exel and A. Vershik}
  {C*-algebras of irreversible dynamical systems}
  {{\it Canadian Mathematical Journal}, to appear, [arXiv:math.OA/0203185]}

\bibitem{\Fowler}
  {N. Fowler}
  {Discrete product systems of Hilbert bimodules}
  {{\it Pacific J. Math.}, {\bf 204} (2002), 335--­375}

\bibitem{\WataTwo}
  {T. Kajiwara and Y. Watatani}
  {C*-algebras associated with complex dynamical systems}
  {preprint, 2003, [arXiv:math.OA/0309293]}

\bibitem{\Kasparov}
  {G. G. Kasparov}
  {Hilbert $C^*$-modules: Theorems of Stinespring and Voiculescu}
  {{\it J. Oper. Theory}, {\bf 4} (1980), 133--150}

\bibitem{\Katsura}
  {T. Katsura}
  {A construction of C*-algebras from C*-correspondences}
  {In {\it Advances in quantum dynamics (South Hadley, MA, 2002)},
Contemp. Math., Amer. Math. Soc., Providence, RI, 2003, pp.~173--182,
[arXiv:math.OA/0309059]}

\bibitem{\Muhly}
  {P. S. Muhly and B. Solel}
  {Tensor algebras over C*-correspondences: representations, dilations, and C*-envelopes}
  {{\it J. Funct. Anal.} (1998), 389--457}

\bibitem{\Pimsner}
  {M. V. Pimsner}
  {A class of C*-algebras generalizing both Cuntz-Krieger algebras and
crossed products by $\bf Z$}
  {{\it Fields Inst. Commun.}, {\bf 12} (1997), 189--212}

\bibitem{\Rieffel}
  {M. A. Rieffel}
  {Induced representations of $C^*$-algebras}
  {{\it Adv. Math.}, {\bf 13} (1974), 176--257}

\bibitem{\RieffelMorita}
  {M. A. Rieffel}
  {Morita equivalence for operator algebras}
  {Operator Algebras Appl., Proc. Symp. Pure Math. vol 38,
pp. 285--298, R. V. Kadison, ed., Amer. Math. Soc., Providence, 1982}

\bibitem{\Takesaki}
  {M. Takesaki}
  {Theory of Operator Algebras I}
  {Springer-Verlag, 1979}

\bibitem{\Watatani}
  {Y. Watatani}
  {Index for C*-subalgebras}
  {{\it Mem. Am. Math. Soc.}, {\bf 424} (1990), 117 pp}

\bibitem{\Zettl}
  {H. Zettl}
  {A characterization of ternary rings of operators}
  {{\it Adv. Math.}, {\bf 48} (1983), 117--143}
